\documentclass[10pt,twoside]{article}

\usepackage{mathtools}
\usepackage{amsfonts} %MIORTEGA
\usepackage{units, amsmath, color}
\usepackage{pstricks,pst-node,pst-coil,pst-plot}  
\usepackage[all,cmtip]{xy}

\usepackage{graphicx}

\usepackage{amssymb}
\usepackage{latexsym}
\usepackage{hyperref}

\makeatletter

\@addtoreset{equation}{section} \makeatother
\newtheorem{theorem}{Theorem}[section]
\newtheorem{remark}[theorem]{Remark}
\newtheorem{lemma}[theorem]{Lemma}
\newtheorem{proposition}[theorem]{Proposition}
\newtheorem{corollary}[theorem]{Corollary}
\newtheorem{definition}[theorem]{Definition}
\newtheorem{example}[theorem]{Example}

\def\1{\mathbf{1}}
\def\:{\lrcorner}
\def\#{\sharp}

\def\a{\alpha}
\def\b{\beta}
\def\g{\gamma}
\def\d{\delta}
\def\e{\varepsilon}

\def\o{\circ}
\def\s{\sigma}

\def\ov{\overline}

\def\qed{\ensuremath{\quad\Box\quad}}

\def\inv#1{\raise.1em\hbox to 0pt{$^{-1}$\hss}_{#1}\;}

\def\V{\noindent}
\def\v{\noindent}

\newcommand{\bean}{\begin{eqnarray*}}
\newcommand{\eean}{\end{eqnarray*}}
\newcommand{\benu}{\begin{enumerate}}
\newcommand{\eenu}{\end{enumerate}}
\newcommand{\eea}{\end{eqnarray}}
\newcommand{\bea}{\begin{eqnarray}}

\def \beit{\begin{itemize}}
\def \eeit{\end{itemize}}
\def \bean{\begin{eqnarray*}}
\def \eean{\end{eqnarray*}}

\def \bui#1#2{\mathrel{\mathop{\kern 0pt#1}\limits^{#2}}}
\def \buil#1#2{\mathrel{\mathop{\kern 0pt#1}\limits_{#2}}}

\newtheorem{Theorem}{Theorem}

\newtheorem{Lemma}{Lemma}

\newcommand{\be}{\begin{equation}}
\newcommand{\ee}{\end{equation}}

\newcommand{\N}{{\mathbb N}}

\newcommand{\R}{{\mathbb R}}
\newcommand{\Q}{{\mathbb Q}}

\newcommand{\SSS}{{\mathbb S}}
\newcommand{\vol}{{\rm vol}}
\newcommand{\Int}{{\rm int}}
\newcommand{\cl}{{\rm cl}}
\newcommand{\pr}{{\rm pr}}

\newcommand{\ben}{\begin{enumerate}}
\newcommand{\een}{\end{enumerate}}
\newcommand{\bit}{\begin{itemize}}
\newcommand{\eit}{\end{itemize}}
\newcommand{\edoc}{\end{document}}

\newcommand{\bdefi}{\begin{definition}}
\newcommand{\btheo}{\begin{theorem}}
\newcommand{\bprop}{\begin{proposition}}
\newcommand{\brema}{\begin{remark}}
\newcommand{\bcoro}{\begin{corollary}}
\newcommand{\blemm}{\begin{lemma}}
\newcommand{\bexam}{\begin{example}}

\newcommand{\edefi}{\end{definition}}
\newcommand{\etheo}{\end{theorem}}
\newcommand{\eprop}{\end{proposition}}
\newcommand{\erema}{\end{remark}}
\newcommand{\ecoro}{\end{corollary}}
\newcommand{\elemm}{\end{lemma}}
\newcommand{\eexam}{\end{example}}

\title{Topologies on the future causal completion}

\begin{document}

\author{Olaf M\"uller\footnote{Institut f\"ur Mathematik, Humboldt-Universit\"at zu Berlin, Unter den Linden 6, D-10099 Berlin, \texttt{Email: mullerol@math.hu-berlin.de}}}

\date{\today}
\maketitle

\begin{abstract}
\v On the Geroch-Kronheimer-Penrose future completion $IP(X)$ of a spacetime $X$, there are two frequently used topologies. We systematically examine $\tau_+$, the stronger (metrizable) of them, which is the coarsest causally continuous topology, obtaining a variety of novel results, among them a complete characterization of the difference in convergence between both topologies. In our framework, we can allow for $X$ being a chr. space and consequently for the interpretation of $IP$ as an idempotent functor on a category that includes spacetimes of very low regularity. Furthermore, we explicitly calculate $(IP(X), \tau_+)$ for multiply warped chronological spaces.
\end{abstract}

\begin{small}{\it Mathematics Subject Classification} (2010):  53A30, 53C50, 53C80\\
\end{small}

\section{Introduction}

\v A central notion of mathematical relativity, frequently used to define concepts related to black holes, is the one of future null infinity. Its definition requires a conformal equivalence of (a part of) the spacetime $M$ to a manifold-with-boundary. Closely related is the recently developed notion of 'conformal future-compact extension' (called 'conformal extension' in some references), which turned out to be useful e.g. for the proof of global existence of conformally equivariant PDEs  \cite{nGoM} like Yang-Mills-Higgs-Dirac equations (the type of the standard model of particle physics), for the proof of existence of black holes in a strong sense in Einstein-Maxwell theory \cite{oM-BH}, as well as in quantization \cite{Dappiaggi}. Let us review the definition following \cite{oM14} and \cite{nGoM}: A subset $A$ of a spacetime $M$ is called {\bf future compact} iff $A $ is contained in $J^-(C)$ for some compact $C \subset M$, and it is called {\bf causally convex} iff there is no causal curve in $M$ leaving and re-entering $A$. A {\bf $C^k$ conformal future-compact extension (CFE)} $E$ of a globally hyperbolic spacetime $(M,g)$ is an open conformal embedding of $(I^+(S_0),g)$ (where $S_0$ is a Cauchy surface of $(M,g)$) into a g.h. spacetime $(N,h)$ with a Lorentzian metric of class $C^k$ such that the closure of the image is future compact and causally convex. This generalizes the usual notion of 'conformal compact extension' (CCE) by requiring only {\em future} compactness of the closure of the image instead of compactness.

\v From the work of Friedrich, Anderson-Chrusciel, Lindblad-Rodnianski and others (\cite{hF}, \cite{AC}, \cite{D}, \cite{DS}, \cite{LR}, \cite{SB}) it follows that there is a weighted Sobolev neighborhood $U$ around zero in the set of vacuum Einstein initial values such that, for any $u \in U$, the maximal vacuum Cauchy development of $u$ admits a smooth CFE (as opposed to CCEs of class $C^2$, whose nonexistence for {\em every} nonflat Einstein-Maxwell solution follows from the positive mass theorem, see Thm. 6 in \cite{oM-BH}). Given a CFE $E: M \rightarrow N$, its future boundary is $\partial^+ ( E) := \{ x \in \partial (E(M)) \vert I^- (x) \cap E(M) \neq \emptyset \} \subset \partial(E(M),N)$. Still, many important examples of spacetimes do not admit CFEs. Luckily, there is a classical {\em intrinsic} notion of future completion $IP(M)$ and future boundary $\partial^+ M$ of a spacetime $M$ as defined by Geroch, Kronheimer and Penrose in \cite{GKP}. Budic and Sachs \cite{rBrS} then defined an appropriate chronological structure $\ll_{BS}$ on $IP(M)$. In Sec. \ref{IP} we give the precise definition $IP(M)$ and $\ll_{BS}$ and revise the well-known construction of a map $i_M: M \rightarrow IP(M)$ assigning to $p \in M$ the set $I^-(p)$, and recall that $IP(M)$ is the set of pasts $I^-(c(\N))$ of {\bf chr. chains} $c: \N\rightarrow M$, that is, a map with $c(m) \ll c(n) \forall m >n$.

Each CFE $E: M \rightarrow N$ yields the {\bf end-point map} $\e_E IP(M) \rightarrow \cl(E(M), N)$ assigning to an IP $I^-(c)$ the future endpoint of $E \o c$ in $N$, defined as the continuous extension of a chr. curve $k : [0; 1) \rightarrow N$ with $k(1-1/n) = E(c(n))$. The choice of $c$ and $k$ is not canonical, but $\e_E$ is well-defined: An endpoint exists by $N$ future-compact, and for two chr. chains $c,k$ with $ I_M^-(c) = I_M^-(k) $ we have $E(M) \cap I_N^-(E \o k) = E(I_M^- (k ) )= E(I_M^- (c)) = I_N^- (E \o c)$, and thus $\e_E(c) = \e_E(k)$ as $(N,h)$ is distinguishing.

\medskip

The present article focuses firstly on {\em future} completions, in contrast to e.g. \cite{FHS1}, and secondly, on the upper part of the causal ladder (assuming at least causal continuity), in contrast to e.g. \cite{H0}. All our constructions can almost effortlessly be made for spaces defined by Harris \cite{Harris1997} called 'chronological spaces' that generalize strongly causal spacetimes and could play a role of limit spaces in Lorentzian geometry similar to metric spaces in Gromov-Hausdorff theory for Riemannian geometry,  cf \cite{oM-Finiteness} for related ideas. Following a definition by Harris, a {\bf chronological} or {\bf chr. set} is a tuple $(X, \ll)$ where $X$ is a set and $ \ll$ is a binary transitive anti-reflexive\footnote{i.e., $\ll$ is disjoint from the diagonal: $\lnot ( x \ll x ) \ \forall x \in X$.} relation on $X$ such that $(X,\ll)$ is {\bf chronologically separable}, meaning that there is a countable subset $S \subset X$ that is {\bf chronologically separating}, i.e., $\forall (x,y) \in \ll \exists s \in S: x \ll s \ll y$. \footnote{Harris' original definition also requires that $\ll$ is {\bf connex}, i.e. $\forall x \in X \exists y \in X: x \ll y \lor y \ll x $. Here we will renounce this requirement as we want to include e.g. compact causal diamonds, which contain points with empty timelike future and past.}

We define $ I^{-} :  X \rightarrow 2^X, I^+(x) := \{ y \in X \vert y \ll x \}, I^-(x) := \{ y \in X \vert x \ll y \} \ \forall x \in X$. For $A \in \{ I, I^+, I^- \}$, we call a chr. set $X$ {\bf $A$-distinguishing} iff $i_X : p \mapsto A(x)$ is injective.

 For a chr. set $X$, a subset $A$ of $X$ is called {\bf past} iff $I^-(A) =A $. A nonempty past subset $U$ of a chr. set $X$ is called {\bf indecomposable} iff for any two past sets $A,B$ we have $U= A \cup B \Rightarrow U =A \lor U = B$, and $X$ is called {\bf preregular} iff $\forall p \in X: I^\pm (p) $ are indecomposable. We call a chr. set $(X, \ll) $ {\bf past-reflecting} iff $\forall x,y \in X: (I^-(x) \subset I^-(y) \Rightarrow I^+(y) \subset I^+(x) )$. If $X$ is a spacetime, past subsets of $X$ have a Lipschitz boundary \cite{BEE}(Th. 3.9), and the purely chr. property of being distinguishing and past-reflecting is equivalent to the mixed chr.-topological property of outer continuity of $I^-$, cf. \cite{eMmS}, Lemmas 3.42 and 3.46. We call $x \in X$ {\bf past-full} resp.\ {\bf future-full} iff $I^-(x) \neq \emptyset $ resp. $I^+ (x) \neq \emptyset$ and {\bf full} iff $I^- (x) \neq \emptyset \neq I ^+ (x) $. We define the category $\mathbf{DPC}$ of $I^-$-distinguishing preregular past-reflecting chr. sets, whose morphisms are $f: M \rightarrow N$ s.t. for $x_1, x_2 \in M$, we have $ f(x_1 ) \ll f(x_2)  \Leftrightarrow x_1 \ll x_2$, and isomorphisms are invertible such\footnote{This is a stronger condition than the notion of isocausality issued by Garc\'ia-Parrado and Senovilla \cite{Senovilla}, inasfar as it requires that the images of the chronological future cones are not subsets of but {\em coincide exactly with} the chronological future cones in the range. However, in \cite{FHS4} it has been shown that the future boundary is not invariant under isomorphisms in the weaker notion from \cite{Senovilla}. To get uniqueness, one can either add the requirement of future continuity of the morphisms  and appeal to Harris' result \cite{Harris1997} or proceed as we do here.}. 

\V It will turn out that we can define a future-compact {\em and} Hausdorff (even metrizable) topology on a future completion if we consider, instead of $IP(X)$, the space $P(X)$ of nonempty past subsets of $X$. The strategy will be to define the causal completion as the closure of the set of pasts of points in $X$ within $P(X)$, and this closure contains, in certain examples, more that just $IP(X)$. 

\V On the set $N(X) := 2^X \setminus \{ \emptyset \}$ of nonempty subsets of $X$ the Budic-Sachs relation $ \ll_{BS}$ with

$$ \forall A,B \subset X: A \ll_{BS} B : \Leftrightarrow \exists x \in B: I^-(x) \supset A . $$

fails to be chronologically separable in general: eg. in the case $X := \R^{1,1}$, we have $ I^-(\{ 0\}) \ll_{BS}  \{ 0\}$, but there is no subset $A$ of $X$ with $ I^-(\{ 0\}) \ll_{BS} A \ll_{BS} \{ 0\}$. It is easy to see that the restriction of $\ll_{BS}$ to the set $P(X)$ of nonempty past subsets {\em is} a chr. relation, and as well the restriction to the set $IP(X)$ of indecomposable subsets of $X$, which we call $\ll_{IP(X)}$. However, it will turn out to be more practical to consider another chr. relation $\ll_+ \supsetneq \ll_{BS}$ on $P(X)$ equally restricting to $\ll_{IP(X)}$ on $IP(X)$, defined by 
 
$$\forall A,B \in P(X) : A \ll_+ B  : \Leftrightarrow \big( \forall U \subset A: U \in IP(X) \Rightarrow U \ll_{BS} B \big) , $$

more on this in the last section. The reason of the prefeence is the following: Whereas in $X:= (-1;1) \times \R \subset \R^{1,1}:= ( \R^2, g:= - dx_0^2 + dx_1^2 )$, $v \ll w :\Leftrightarrow g(v-w, v-w ) <0 \land v_0 < w_0$, the past subset $A_a := I^- (\{ 0\} \times (-a; a))$ has empty chr. future w.r.t. the chr. relation $\ll_{BS}$ for all $a >1$, the future w.r.t. $\ll_+$ contains e.g. the sets $ I^- (\{ r\} \times (-a;a))$ for all $r >0$. Consequently, if we later defined the topology $\tau_+$ on $P(X)$ by means of the chr. relation $\ll_BS$, we would not get convergence of $A_{a+1/n}$ to $A_a$ for $n \rightarrow \infty $ as desired, whereas we do if we work with the chr. relation $\ll_+$ on $P(X)$.

For every object $X$ of $C$ we have the morphism 

\[ i_X: X \rightarrow IP(X), x \mapsto I^- (x) \forall x \in X \] 

which we can enrich by chr. relations to define a functors \\ $ P_+: X \mapsto (P(X), \ll_+)  $ and $IP : X \mapsto (IP(X), \ll_{IP(X)})$ from $C $ to $ C$.  
 
\medskip

There are functorial ways to define a causal relation from a chr. relation $\ll$, e.g. $\a$ defined by 

\[ x \ \a (\ll) \   y : \Leftrightarrow I^+(y) \subset I^+(x) \land I^-(x) \subset I^-(y) , \]

 see \cite{eMmS}, Def. 2.22 and Th. 3.69 showing that for a causally simple spacetime $(M,g)$, with the usual definition of $\ll$ resp. $\leq$ via temporal resp. causal curves and with the definition $J^+(x) := \{ y \in X : y \geq x\}$ and $J^-(x) := \{ y \in X: y \leq x \} $, we get $x \leq y \Leftrightarrow x \a (\ll) y \Leftrightarrow I^+(y) \subset I^+(x) \lor I^-(x) \subset I^-(y)$: The implication from left to right is true in any spacetime, for the other direction calculate $ I^+(y) \subset I^+(x) \Rightarrow y \in \cl (I^+(x))  = \cl (J^+(x) ) = J^+(x)$. It is well-known that $\alpha$ may introduce spurious causal relations in the non-causally simple case, see the causally continuous example $\R^{1,n} \setminus \{ 0 \}$, where $(-1,-1) \not\leq (1,1) $ but $(-1,-1) \a (\ll) (1,1)$. We will also see in Sec. \ref{posets} that the definition of the causal relation via $\a$ entails the {\em push-up property} $I^+ \o J^+ = I^+$.

 Thus, until Sec. \ref{posets} we just renounce speaking of causality and restrict ourselves to the chr. relation.

Let $Y$ be a chr. set. A map $c: \N \rightarrow Y$ is called (future) chr. {\bf chain} iff, for all $m, n \in \N$, we have  $m < n \Rightarrow c(n) \ll c(m)$. A chr. chain $c$ has a {\bf (chr.) limit} $x \in X$ iff $I^-(x) = \bigcup \{ I^- (c(n)) \vert n \in \N \}$. Analogously to chr.\ chains, (future) chr.\ curves are defined by replacing $\N$ with a real interval.  $Y$ is called {\bf future properly causal} (cf. \cite{FHS1}) or {\bf future chr. complete} iff for every future chr. chain $n \mapsto c_n $ in $Y$ there is $c_\infty \in Y $ with 
$I^-(c_{\infty}) = \bigcup_{i=0}^{\infty} I^{-} (c_i)$, and it is a classical fact that $IP(X)$ is future chr. complete.  If $X$ is preregular, the map $i_X: p \mapsto I^-(p)$ is a morphism of $C$ (due to chr. separability) with $i_X(X) = I^-(\partial^+ (IP(X))$, where $\partial^+ (Y) := \{ y \in Y \vert I^+(y) = \emptyset \}$ for a chr. set $Y$.  As we show in Sec. 2, {\em for every preregular chr. set $X$, the set $(IP(X), \ll_{IP(X)})$ is  past-reflecting}. 

\medskip

A 3-tuple $(X, \tau, \ll )$ where $(X, \tau)$ is a topological space and $ (X, \ll) \in {\rm Obj} (C)$ is called a {\bf chronological} or {\bf chr. space}. It is called {\bf regular} iff it is preregular and $ \forall p \in X: I^\pm (p)  $ open. In any preregular chr. space, past subsets are open. Causal spacetimes are regular chr.\ spaces. We define the category $CS$ of distinguishing, past-reflecting and regular chr. spaces whose morphisms are the continuous morphisms of $C$ and whose isomorphisms are the invertible such.
We enrich the functor $P_+: X \mapsto P(X)$ by adding one of two topologies, i.e., prolong $P_+$ to a functor $F_\pm = \tau_+ \o IP$ with a functor $\tau_\pm: C \rightarrow CS$ which is a left inverse of the forgetful functor $CS \rightarrow C$. We show that both choices $F_\pm$ satisfy three desiderata: Both are {\em future chain complete}, {em marginal} and {\em respect CFEs}, where, for a functor $F: C \rightarrow C$:

\begin{itemize}
\item $F$ is {\bf future chain complete} iff for each future chr. chain $c$ in each $X  \in {\rm Obj} (F_\pm(C))$, $c(n) \rightarrow_{n \rightarrow \infty} c_\infty$;
\item $F$ is {\bf marginal} iff for every causally continuous spacetime $M$. the map $i_M$ is a homeomorphism onto $i_M (M) $, which is open and dense in $F(M)$;
\item $F$ {\bf respects CFEs} iff the end-point map $\e_E$ of any CFE $E: M \rightarrow N$ of $M$ maps $\partial^F(M) := F(M) \setminus i_M(M) $ homeomorphically to $\partial^+ (E)$, preserving causal relations.
\end{itemize}

The history of topologies on $IP(M)$ is strangely involved. The first proposal for a topology, here denoted by $\tau_+$, was by Beem \cite{Beem!}, for the entire (future and past) causal completion of a causally continuous spacetime, which consists of the union $IP(X)$ and its time dual modulo some identifications  --- the proofs in Beem's article however transfer verbatim to $IP(X)$. He stated without  proof that $\tau_+$ is compatible with the metric defined in Th. \ref{TFAE}, Item 2, and showed:

\begin{enumerate}
	\item If $X$ is a causally continuous spacetime then $\tau_+$ is metrizable on $IP(X)$ and $i_X(X)$ is open and dense in $IP(X)$;
	\item If $X$ is g.h., then $i_X^{-1} (\tau_+)$ is the Alexandrov (manifold) topology on $X$. 
\end{enumerate}

After Beem's article, $\tau_+$ seems to be have fallen into oblivion in the Lorentzian community, considering that in the 1990s a lively and extensive  debate about the choice of topology on $X$ began without ever taking $\tau_+$ into account. Another very agreeable candidate, the 'chronological topology', here denoted by $\tau_-$, was elaborated by Flores \cite{jF} inspired by ideas of Harris \cite{H0}. Flores, Herrera and S\' anchez  \cite{FHS1} extended its definition to the entire causal completion and related it to the conformally standard stationary situation \cite{FHS2}. $\tau_-$ is in general not Hausdorff --- but still T1. The main advantage of $\tau_-$ is that it recovers the manifold topology on a spacetime $M$ not only if $M$ is causally continuous (like $\tau_+$ does), but in the larger class of strongly causal spacetimes.

In 2014, unaware of Beem's article, the author defined a metric on $IP(X)$ as in Th. \ref{TFAE}, Item 2, using it as a tool in the context of conformal future-compact extensions to show that the future boundary of any CFE is homeomorphic to the intrinsic future boundary $\partial^+ X$ defined in terms of $IP(X)$. 
Then in 2018 Costa e Silva, Flores and Herrera, also unaware of Beem's result, explored $\tau_+$, re-proving many facts independently and more systematically, and contributed other interesting aspects, e.g. the fact that if $(IP(X), \tau_-)$ is Hausdorff then $\tau_- = \tau_+$, and applied the causal boundary topologized by $\tau_-$ to generalize the notion of black hole showing a version of Hawking's theorem that each trapped surface is contained in the complement of the past of future null infinity under appropriate conditions on the future boundary. 
Shortly after publication of their article, the author and they discussed the topology again discovering Beem's article. This led the author to work out a streamlined introduction into the topic, comparing the topologies systematically. While writing this review article, the author obtained several novel results, also laid down here.

Whereas many previous approaches defined $\tau_+$ only on $IP(X)$, we present a definition of $\tau_+$ for each chr. set., give a self-contained introduction to $IP(X)$, compare $\tau_+$ and $\tau_-$, showing, e.g.:

\begin{enumerate}
\item several steps of an analogue of the causal ladder for chr. spaces (cf. Sec. \ref{ladder} and Sec. \ref{posets});
\item the fact that for all causally continuous {\em well-behaved} \footnote{defined below, satisfied e.g. by each connected causal spacetime and each connected Lorentzian length space} chr. space $X$, the chr. space $\cl (i_X(X), (P(X), \tau_+)) \supset IP(X)$ is well-behaved and causally continuous again (Th. \ref{recovering}), 

\item equivalence of several definitions of $\tau_+$, some of which appear in the previous literature, among them the equivalence stated without proof in Beem's article (Th. \ref{TFAE});
\item future-compactness of $\tau_+$ on $P(X)$ and implications like existence of minimal and maximal TIPs in such futures (Th. \ref{ConvergencePlusMinus}).
\end{enumerate}

Moreover, as an example we calculate $(IP(X), \tau_+)$ for $X$ being a {\em multiply-warped chr. space} (see below, generalizing multiply-warped spacetimes), extending a result by Harris (see Sec. \ref{multiply-warped}).

\medskip

The article is structured as follows:  Sections 2, 3, 4 are largely expository; Section 2 reviews the construction of the functor $IP$ assigning to a chronological set $M$ its future completion $IP(M)$ as a chr. set, Section 3 gives an overview over limsup and liminf for subsets of ordered spaces. Section 4 resumes facts about $\tau_-$ on $IP(M)$. Section 5 explores parts of the causal ladder for chronological spaces. Section 6 introduces the topology $\tau_+$ on a chr. set, compares it to $\tau_-$ and to the Alexandrov topology and applies it to $IP(M)$. In Section 7, in a close analogy to a proof by Harris, we compute $\tau_+$ in the future completions of multiply warped products. In Section 7, we develop the topic from the perspective of the causal relation instead of the chr. relation, develop a functor $\tilde{\a} \neq \a$ from the chr. category to a causal one, show that $\tilde{\a}$, when applied to arbitrary causally continuous spacetimes, avoids the spurious causal relations known from $\a$, and give a simple one-line definition of the future causal completion as a causal space with topology $\tau_+$, and finally, draw a a conclusion passing a solomonic judgement between the two topologies.

\medskip

\V The author is very grateful for several valuable discussions with Iv\'an Costa e Silva, Jos\'e Luis Flores, Leonardo Garc\'ia-Heveling, Stacey Harris, Jonat\'an Herrera and Miguel S\'anchez, and for very detailed and helpful reports of anonymous referees on previous versions of this article.

\section{Sets of past subsets}
\label{IP}

Let $X$ be a chr. set. Arbitrary unions of past subsets of $X$ are past. $I^-(A)$ is past for each $A \subset X$, i.e., $I^-$ is idempotent: Transitivity of $\ll$ implies $I^-(I^-(A)) \subset I^-(A)$ and chr. separability implies the reverse inclusion. Geroch, Kronheimer, Penrose \cite{GKP}(Th. 2.1) and Harris \cite{Harris1997} show:

\medskip

{\em A nonempty subset of a chr.\ set $X$ is an IP if and only if it is the past of the image of a chr.\@ chain.} 

\medskip

For the sake of self-containedness let us give a proof: For the 'if' part, let $ B \cup C = A = I^- (c(\N))$ for a chr. chain $c$ and $B,C  $ past. As $c^{-1} (B) \cup c^ {-1} (C) = \N  $, one of them is infinite, and this one equals $A$. For the other direction, call a chr. space $Y$ {\bf (past) synoptic} iff for all $p,q \in Y$ we have $I_Y^+(p) \cap I_Y^+(q) \neq \emptyset$. Then we prove that {\em a past set is synoptic if and only if it is indecomposable}: Nonsynopticity implies existence of two points $p,q$ with $I_Y^+(p) \cap I_Y^+(q) = \emptyset$, or equivalently, $q \notin U_Y^+(p):=  I_Y^-(I_Y^+(p)) $, and a nontrivial decomposition of $Y$ is then $(U_Y^+(p), X \setminus I_Y^+(p))$; for the other direction let a nontrivial decomposition $(A,B)$ be given, then for $p \in Y \setminus A, q \in Y\setminus B$ one easily verifies $I_Y^+(p) \cap I_Y^+(q) = (I_Y^+(p) \cap I_Y^+(q)) \cap (A \cup B) = \emptyset$. Knowing synopticity, it is easy to construct a chr. chain $q$ generating the IP: Let $p(\N)$ be chr.ly separating in $Y$, then choose inductively $q_{i+1}$ from $  I^+(p_i ) \cap I^+(q_i) $ (nonempty, as $Y$ is synoptic). Then for all $x \in Y$ we know, again by separability, that $I^+(p) $ contains some $p_i \ll q_{i+1}$, and so $Y= I^-( q(\N))$.

If $X$ is preregular, then also finite intersections $ \bigcap_{k=1}^n A_n$ of past subsets $A_n$ are past: By induction it suffices to consider $n=2$, and then for $p \in A_1 \cap A_2$ there are $p_i^+ \in I^+ (p) \cap A_i$, and future synopticity of $I^+ (p)$ implies that we can find $p^+ \in I^+(p) \cap I^- (p_1^+) \cap I^-(p_2^+) \subset I^+(p) \cap A_1 \cap A_2$. On the other hand, in each spacetime there are countable intersections of past subsets that are not past: For $A_n := I^-((0, 1/n)) \subset \R^{1,1}$, we get $\cap_{k=1}^\infty A_n = I^-(0) \cup [0; \infty) \cdot (-1,1)$.

We define, for any nonempty subset $R   $ of $X$, by $I^{+ \bigcap}_X(R)$ the {\bf joint} or {\bf common chr. future}  $\bigcap \{ I^+(r) \vert r \in R \}$ {\bf of $A$ in $X$}, correspondingly for the pasts. Whereas $I^\pm$ are monotonically increasing maps from $(N(X), \subset) $ to itself, $I^{\pm \bigcap}$ are monotonically {\em decreasing}, and e.g. $I^{+ \bigcap} (A \cap B) \supset I^{+ \bigcap} (A) \cup I^{+ \bigcap} (B)$. We then have $A \ll_{BS} B \Rightarrow A \subset B$ and

 $$A \ll_{BS} B   \Leftrightarrow I^{+  \bigcap}_X (A) \cap B \neq \emptyset \ \forall A,B \in 2^X.$$

Essentially as in \cite{jF} (Th.4.3), we see that $IP(X)$ is $I^-$-distinguishing if $X$ is $I^-$-distinguishing, as $i_X^{-1} (I^-_{P(X)} (x)) = x   \ \forall x \in IP(X)$. Equally $(P(X), \ll_+)$ is $I^-$-distinguishing, even more, {\em the map $I_{P(X)}^- \cap IP(X): P (X) \rightarrow IP(X)$ is injective}: Let $A,B \subset X$ past with $ I^-_{P(X)} (A) \cap IP(X) = I^-_{P(X)} (B) \cap IP(X) $, thus for all $D \in IP(X): I^{+ \cap } (D) \cap A \neq \emptyset \Rightarrow I^{+ \cap} (D) \cap B \neq \emptyset $ (*); now let $x \in A = I^-(A)$, then there is $y \gg x$ with $y \in A$, and if we apply (*) to $D:= I^-(x)$ with $y \in I^{+ \cap} (D) \cap A \neq \emptyset$, there is $b \in B \cap I^{+ \cap } (D) $ and thus $I^-(x) \subset I^- (b) $, and by $X$ past-reflecting $I^+(b) \subset I^+(x)$. Now by $B$ past, there is $c \in B$ with $c \in I^+(b) \subset I^+(x)  $, thus $x \in I^-(c) \subset B$.

{\em If$(X, \ll)$ is an $I^-$-distinguishing preregular chr.\ space, $i_X:(X, \ll)  \rightarrow (IP(X),  \ll_+ )$ is an injective morphism of $\mathbf{DPC}$ (i.e., preserves the chr.\  relation). The functor $IP$ assigning $IP(X)$ to any chr. set $X$ is idempotent (as it maps each future-complete chr. space to itself), it maps (preregular resp.\ past-full resp.\  $I^\pm$-distinguishing) chr.\ sets to future complete (preregular resp.\ past-full resp. $I^\pm$-distinguishing) chr. sets, $i_X $ preserves limits of chr.\ chains, and $i_X(X)$ is past in $IP(X)$.}

This is shown via $ \forall A \in P(X): A = \bigcup \{ C \in IP(X) \vert C \ll_+ A \} $ to show past-distinguishing; past-reflecting follows via $I^-_{BS} (A) \subset I^-_{BS} (B) \Leftrightarrow A \subset B  $, shown in turn as follows: Assume $\emptyset \neq A \setminus B \ni p $, then there is $p_+ \in I^+(x) \cap (A \setminus B)$; by $C:= I^-(p)$ get $ C \in I^-_{BS}( A) \setminus I^-_{BS} (B)$, contradiction.

\newpage

\section{liminf and limsup on posets}

In this article we use liminf and limsup in several instances. To avoid confusion about those, this section serves to recall that all of these definitions actually are the same. Let $(Q, \leq) $ be an ordered set and let $A \subset Q$. Throughout this article, let $x<y : \Rightarrow x\leq y \land x \neq y$. Recall that $q \in Q$ is called {\bf upper} resp. {\bf lower bound of $A$} iff $a \leq q \forall a \in A$ resp. $a \geq q \forall q \in Q$, and $q$ is called {\bf supremum} resp. {\bf infimum of $A$} if $q$ is an upper resp. lower bound of $A$ and iff there is no upper resp. lower bound $\check{q}$ resp. $\hat{q}$ of $A$ with $\check {q} < q$ resp. $\hat{q} > q$. It is well-known and easy to show that infima and suprema of subsets $A$ of posets, whereas they do not always exist, are always unique. Thus $\inf$ and $\sup$ define partial maps from $P(Q) $ to $ Q$. For each sequence\footnote{We can even include nets by changing $\N$ for some directed set $N$.} $a: \N \rightarrow Y$,

\bea
\label{Limsupinf} \liminf (a) := \sup_{n \in \N} \inf_{m \geq n} a(m), \qquad \qquad  \limsup (a) := \inf_{n \in \N} \sup_{m \geq n} a(m) ,
\eea

if the involved suprema and infima exist. Thus $\limsup$ and $\liminf$ are partial maps from $X^\N$ to $X$. This will be the only definition of $\liminf$ and $\limsup$ we use in the article. Writing $M \nearrow N$ for a strictly increasing map between posets $M$ and $N$, we now specialize the definition to 3 posets:

\bigskip

{\bf 1.} In the case of subsets $Y= 2^X$ of a set $X$ ordered by inclusion, we get $\inf = \bigcap$, $\sup = \bigcup$ and

\bean
\liminf_{2^X}(a) &:=& \bigcup_{n \in \N} \bigcap_{m \geq n} a(m), \\
\limsup_{2^X}(a) &:=& \bigcap_{n \in \N} \bigcup_{m \geq n} a(m) = \bigcup \{ \bigcap_{n \in \N} (a \o j) (n) \  \big\vert \  j: \N \nearrow \N   \}.
\eean 

{\bf 2.} Applied to the case $C(X)$ of closed subsets of a topological space $X$, ordered by inclusion, we get $\inf (A) = \bigcap A$ and $\sup (A) = \cl (\bigcup A)$. 
 For $A: \N \rightarrow 2^X$, the defining equality \ref{Limsupinf} yields

\bean
\limsup \  _{C(X)} (A) & :=& \{ q \in X \vert  \exists a: \N \rightarrow X \exists j: \N \nearrow \N : a \o j (n) \in A(j(n)) \forall n \in \N \land a(j(n)) \rightarrow_{n \rightarrow \infty} q \}, \\
\liminf \ _{C(X)} A  &:=& \{ q \in X \vert  \exists a: \N \rightarrow X : \forall n \in \N\exists k \in \N \forall m \geq k : a(n) \in A(m) \land a(n) \rightarrow_{n \rightarrow \infty} q \}.
\eean

For $a: \N\rightarrow X^\N$ we define $\widetilde{\liminf} (a) := \liminf (\cl \o a) \in C(X)$ and $\widetilde{\limsup} (a) := \limsup (\cl \o a) \in C(X) $. With those definitions, we see that 

$$\widetilde{\liminf} (a) \supset \cl (\liminf a ) , \  \widetilde{\limsup} (a) \supset \cl (\limsup a) .$$ 

\V Even for sequences of open subsets, these inclusions are not equalities: Let $a(n) := (\frac{1}{n+1}; \frac{1}{n+2}) \subset \R$, then the $a(n)$ are pairwise disjoint, so $\limsup(a)  = \emptyset$ but $\widetilde{\limsup} a = \{ 0 \}$. For more details see \cite{Busemann}.

\medskip

{\bf 3.} On the set $P(X)$ of past sets of a chr. set ordered by inclusion, $\sup (A) = \bigcup A$, $\inf (A)  = I^-(\bigcap A)$,

\[ \liminf (a) = \bigcup_{n \in \N} I^-(\bigcap_{m \geq n} a(m)) = I^- ( \overbrace{\bigcup_{n \in \N} \bigcap_{m \geq n} a(m)}^{= \liminf_{2^X} a }) , \  \ \limsup (a)  =  I^-(\limsup \ _{2^X} a) . \]

\newpage

\section{The topology $\tau_-$}

We define a topology $\tau_-$ by letting $C \subset IP(M)$ be $\tau_-$-closed iff for every sequence $\s$ in $C$ we have 

$$ C \supset L_-(\s) := \{ P \in IP(M) \setminus \{ \emptyset \} \vert P \subset \liminf (\s) \land P {\rm \ maximal \ IP \ in \ } \limsup (\s) \} .$$

where $\limsup$ and $\liminf$ are defined set-theoretically, see Sec. \ref{ladder}. There is a construction that unites future and past boundary, and then, one can single out the future part of the boundary and equip it with the relative topology; the latter coincides with the future chronological topology if the initial spacetime $(M,g)$ was globally hyperbolic \cite{FHS1}. The definition of $\tau_-$ works for $IP(X)$, where $X$ is any chronological space, and is functorial in the category of $IP$s of chronological spaces and the IP prolongations $I^- \o f$ of chronological morphisms $f$. Furthermore, $\tau_-$ is locally compact and even future-compact, which follows from the fact shown in \cite{FH}(Theorem 5.11) that every sequence of IPs not converging to $\emptyset$ has a subsequence convergent to some IP. Defining a functor $F_-:= (IP, \tau_-)$ between the categories of causally continuous spacetimes and chr. spaces we get:

\begin{Theorem}[see \cite{FH}, \cite{FHS1}]
$F_-$ is sequentially future-compact, marginal and respects CFEs.
\end{Theorem}

\V {\bf Proof.} Sequential future compactness is shown in \cite{FH}, Th. 5.11, marginality in  \cite{FHS1}, Th. 3.27. Now, Theorem 4.16 in \cite{FHS1} assures that the end-point map is a chronological homeomorphism if 

\begin{enumerate}
\item $E$ is future chr. chain complete,
\item each point $p \in \partial^+ E(M) $ is {\bf timelike transitive}, i.e. there is an open neighborhood $U$ of $p$ in $N$ such that for all $x,y \in V:= \cl (E(M)) \cap U$ the following push-up properties hold: $x \ll_E z \leq_E y \Rightarrow x \ll_E y$, $ x \leq_E z \ll_E y \Rightarrow x \ll_E y$. Here, $p \ll_E q$ iff there is a continuous curve $c$ in $\cl (E(M))$ between $p$ and $q$ that is a smooth future timelike curve in $M$ apart from the endpoints of the interval.    
\item Each point $p \in \partial^+ E(M) $ is {\bf timelike deformable}, i.e. there is an open neighborhood $U$ of $p$ in $N$ such that $I^-(p, E(M)) = I^-(E(c))$ for a $C^0$-inextendible future timelike curve $c: [a;b] \rightarrow M$ with $E(c(a)) \ll_E q $ for all $q \in V:= \cl (E(M)) \cap U$. 
\end{enumerate} 

These properties are satisfied by CFEs: The first one follows from future compactness, the second one from the push-up properties in $N$ and the fact that $\leq_{E(M)} \subset \leq_N$ and for $q \in \partial^+ E(M) $ we have $p \ll_{E(M)} q  \Leftrightarrow p \ll_N q$ as $I^-(q) \subset E(M)$. The same argument works also for the last item. \hfill \qed

\bigskip

\V The topology $\tau_-$ does not in general inherit the $\R$-action of the flow of the timelike Killing vector field even for standard static spacetimes, without the additional hypothesis that $\tau_-$ is Hausdorff.    
Remark 3.40 in \cite{FHS1} shows that $\tau_-$ is in general not Hausdorff  nor even first countable; in Sec. \ref{MetricsOnIPs} we will see an example (the ultrastatic spacetime over the unwrapped-grapefruit-on-a-stick) well-known to have a non-Hausdorff future causal boundary if the latter is equipped with $\tau_-$.  

\newpage

\section{Topologies on power sets, causal ladder for chr. spaces} 
\label{ladder}
    
As $IP(X)$ consists of subsets of $X$, we now examine more closely topologies on power sets. For topological spaces $(W, \s),(X, \tau)$, a map $F: W \rightarrow \tau(X)$ from $W$ to the open sets of $X$ is called {\bf inner (resp., outer) continuous at a point $p \in Y$} iff for all compact sets $C \subset (F(p)) $ (resp., for all compact sets $ K \subset X \setminus \cl (F(p))$), there is a $\s$-open set $U$ containing $p$ such that for all $q \in U$, we have $C \subset (F(q))  $ (resp., $K \subset X \setminus \cl  (F(q) )$). Translated to notions of convergence, this means that for every net $a$ valued in $\tau(W)$ convergent to $p$, the net $F \o a$ {\bf inner (resp., outer) converges to $F(p)$}, where a net $a: b \rightarrow \tau (X)$ inner converges to $A \in \tau( X)$ iff for all compact sets $C$ in $A $ there is $n \in b$ such that for all $m \geq n$ we have $C \subset a(m) $ (resp., for all compact sets $C$ in $\Int (X \setminus A) $ there is $n \in b$ such that for all $m \geq n$ we have $C \subset \Int (X \setminus a(m)) $). We denote the topology determined by inner and outer convergence (i.-o.-convergence) by $\tau_{io}$ (which is, somewhat confusingly at a first sight, a topology on the set $\tau(X)$ of open subsets of $X$). Sometimes we will consider the i.o. topology as a topology on the entire power set $2^X$; there, it is induced as the initial topology $\Int^{-1} (\tau_{io})$ via the map $\Int: 2^X \rightarrow \tau (X)$ taking the interior w.r.t. $\tau$. A net $a$ converges in this topology to a subset $U$ iff $ \Int \o a$ converges to $\Int (U)$. Of course, this topology is not Hausdorff any more: For each $A \in \tau (X)$, $A$ and $\cl (A)$ cannot be separated by open neighborhoods.

\V Let $C(X)$ resp. $K(X)$ denote the set of closed resp. compact subsets of $X$; if a $X$ is a metric space, until further notice we always equip $K(X)$ with the topology defined by the Hausdorff distance $d_H$ (recall that the definition of the Hausdorff distance implies that $d_H(\emptyset, A) = \infty$, which makes the Hausdorff distance a generalized metric). Recall that a metric space $X$ is called {\bf Heine-Borel} iff closed bounded subsets are compact. In each Heine-Borel space $X$, each compact $K \subset X$ and each $R>0$, the subset $B(K,R)$ is compact.

\begin{Theorem}
\label{HausdorffInnerOuter}
Let $X$ be a Heine-Borel metric space, let $a: Y \rightarrow  \tau(X)  $ be a net. If $ a$ i.-o.-converges to $a_\infty \in \tau (X)$, then $  \forall K \subset X {\rm \ compact } : a_K := \cl \o (\cdot \cap K)  \o a: Y \rightarrow (K(X), d_H)  $ converges to $\cl (a_\infty \cap K)$.

\end{Theorem}
{\bf Remark.} The reverse implication does not hold, an instructive counterexample with $a$ taking values in the interiors of closed connected subsets being $a(n):= B(0,1) \setminus B( x(n), 1/n)  \subset \R^k$ for $k \geq 2$, where $x: \N \rightarrow \Q^k \cap B(0,1) $ is any bijective sequence; in this case, for $K:= \cl (B(0,1))$, the sequence $a_K  $ converges to $\cl (B(0,1) )$ in $d_H$, but $a$ does not converge in $\tau_{io}$.  

\V {\bf Proof}. We have to show $\cl (K \cap a(  m)) \rightarrow_{m \rightarrow \infty} \cl (K \cap a_\infty)) $ in $d_H$. Inner convergence implies that for all $L \subset a_\infty$ compact, there is $m_L \in Y$ such that $L \subset a(n) \forall n \geq m_L$, and outer convergence implies that  for all $S \subset \Int (X \setminus a_\infty)$ compact, there is $\mu_S \in Y$ such that $S \subset \Int (X \setminus  a(n)) \forall n \geq \mu_S$. Apply this to $L := K \setminus B(X \setminus a_\infty, \epsilon)$ and $S := \ov{B(x, \e/2)} \setminus B(a+\infty, \e/2)$, showing the claim. \hfill \qed

\medskip

We can express many facts about noncompact subsets in a more elegant way via Busemann's following result (see \cite{Busemann}) on a metric $d_1$ on the space $C(X)$ of the closed nonempty subsets of a Heine-Borel metric space $(X,d_0)$, which differs from Hausdorff metrics like the continuous extension of $\arctan \o d_H$ from $d_H$, by blurring the behavior at infinity:

\begin{Theorem}[Compare with \cite{Busemann} Sec. I.3]
\label{Busemann}
Let $(X, d_0) $ be a Heine-Borel metric space, let $x_0 \in X$. Then for $C(X) := \{ A \subset X \vert  \emptyset  \neq A {\rm \ closed }\} $, the map 

\bean
&d_1 = d_1^{x_0, d_0}:& C(X) \times C(X) \rightarrow \R \cup \{ \infty \} ,\\
&d_1(A,B)& := \sup \{  \underbrace{\vert d_0 (  x ,A) - d_0 ( x,B) \vert \cdot \exp (-d_0( x_0 , x))}_{\psi_{A,B} (x)} ; x \in X \}  \ \forall A,B \in C(X)
\eean

takes only finite values and is a metric on $C(X)$ that makes $C(X)$ a Heine-Borel metric space. For every point $y$ we have $d_1 (\{x_0\}, \{ y \}) = d_0 (x_0,y)$. For two points $x_0, x_0'$, the metrics $d_1^{x_0}$ and $d_1^{x_0'}$ are uniformly equivalent (so by abuse of notation we henceforth occasionally suppress this dependence). A sequence $a: \N \rightarrow 2^X$ converges in the topology induced by the extended metric $d_1$ to $a_\infty$ iff $ \widetilde{\liminf}_{n \rightarrow \infty} a(n)  =  \widetilde{\limsup}_{n \rightarrow \infty} a(n) = \cl (a_\infty)$. Limits w.r.t. $d_1$ are monotonous in the sense that if two $d_1$-convergent sequences $a$ and $b$ in $C(X)$ with limits $a_\infty$, $b_\infty$, respectively, satisfy $a(n) \subset b(n) \forall n \in \N $ then $a_\infty \subset b_\infty$ (which entails an obvious sandwich principle for convergence). \\ 
For each $y \in X$ there is $K>0$ s.t. for all $A,B \in C (X)   $ we have $y \in A \cap B  \Rightarrow d_1^{x_0,d_0}(A,B) < K$.
\end{Theorem}

\V {\bf Proof.} The reference \cite{Busemann}, Sec.~I.3 shows all assertions but the last one. With $c:= d_0(x_0,y)$ we get

\bean
\vert d_0 ( x ,A) - d_0(x,B) \vert \leq \max \{ d_0(x, A), d_0(x,B) \} \leq \underbrace{d_0(x,x_0)}_{=:s} + \underbrace{d_0(x_0,y)}_{=:c} .  
\eean

$K:= \sup \{ (s+c) \cdot e^{-s} \vert s \in [0;\infty) \} < \infty$, and $\vert d_0 ( x,A) - d_0( x,B) \vert \cdot e^{-d_0(x_0,x)} < K \ \forall x \in X$. \\ In the calculation of $d_1(\{ x_0\}, \{ y \})$, for '$\geq$' take $x=x_0$, for '$ \leq$' calculate \\ $(d_0 (x,x_0) - d_0 (x,y))e^{-d_0 (x_0, x)} \leq d_0(x_0, y) e^{-d_0(x_0, x)} \leq d_0(x_0, y)$. \hfill \qed

\begin{Theorem}
\label{HB}
Let $(X,d)$ be a Heine-Borel metric space. Then $d_1$-convergence of a sequence $a: \N \rightarrow C(X)$ is equivalent to $d_H$-convergence of $a_L :=    L \cap  a  $ for all $L \subset X$ compact. Thus $f: Y \rightarrow (C(X), d_1)$ is continuous iff $f_L: Y \rightarrow (K(L), d_H) $ (as in Th. \ref{HausdorffInnerOuter}) is continuous for all $L \in K(X)$.
\end{Theorem}

\V {\bf Proof.} The first direction, to show that $d_1$-convergence of $a$ implies $d_H$-convergence of $a_L$, follows from the topological equivalence of $d_1$ and $d_H$ on compacta. For the reverse implication, assume $a_K(m) \rightarrow_{m \rightarrow \infty} b_K$ for some $b \in C(X)$, let $\e >0$, then choose $K:= \cl (B(x_0, r))$ for $2r e^{-r} < \e$ and $n\in \N$ such that for all $m \geq n$ we have $d_H (a_K(m), b_K) < \e$. Then $d_1 (a(m), b) < \e $.   \hfill \qed

\medskip

Now we want to connect the previous notions to causality. For a sequence $a$ of past sets, in general neither $\liminf (a)$ nor $\limsup (a)$ is a past set: Consider a Lorentzian product $X:= (0;1) \times M$, let us take $M:= (0;1)$ for simplicity. As $(0;1) \times (0;1) $ is isometrically  embedded in $\R^{1,1}$, the future causal completion is isomorphic in the category $C$ to the subset $ V:= (0;1] \times [0;1] $ of $\R^{1,1}$. Now consider the sequence $a $ in $V$ defined by $a(n) := (1, \frac{1}{2} + \frac{1}{n})$ and a corresponding sequence $A$ of indecomposable past subsets of $X$ defined by $A(n):= I^-(a(n))$. Then we get $\liminf (A) = I^-((1, \frac{1}{2})) \cup \{ (1, \frac{1}{2}) + t(-1,1) \vert 0< t< \sqrt{2} \}$, and the second subset in the union is a part of the boundary of the first subset,  preventing the liminf to be open and so to be a past set. Analogous examples exist in Kruskal spacetime. However, we can define  a subset $B$ of a chr. set $X$ to be {\bf almost past} resp. {\bf almost future} iff $I^-(B) \subset B$ resp. $I^+(B) \subset B$ (if $X$ is a regular chr. space, these are equivalent to $I^-(B) \subset \Int (B)$ resp. $I^+(B) \subset \Int (B)$). We denote the set of all closed almost past resp. almost future subsets of $X$ by $CAP(X)$ resp. $CAF(X)$.   We call a chr. space $X$ {\bf chr. dense} iff $ p \in \cl (I^\pm(p))  \ \forall p \in  X  \setminus \partial^\pm X  $. Then we get easily:

\begin{Theorem}  
\label{AlmostPast}
Let $X$ be a chr. set and $a: \N\rightarrow P(X)$. 
\begin{enumerate} 
	\item $\limsup_{2^X}a $ and $\liminf_{2^X} a $ are almost past.
	\item If $X$ is regular, then $\widetilde{\limsup} (a) = \cl (\limsup a)$ and $\widetilde{\liminf} a  = \cl (\liminf a)$. 
	\item If $X$ is regular and chr. dense, then the complement of an almost past subset is almost future, each almost past or almost future set $B$ satisfies $B \subset \cl(\Int B)$, thus $ \cl(\Int B) = \cl (B)$, $\Int (B) = \Int(\cl(B)$, and the maps $\cl: P(X) \rightarrow CAP (X)$ and $I^-: CAP (X) \rightarrow P(X)$ are inverse to each other. Moreover, $I^-(\Int(A)) = I^-(A)$ for each $A \in AP(X)$.
	\end{enumerate} 
\end{Theorem}

{\V Proof.} Let $x \in I^- (\liminf_{2^X} (a)) = I^-(\bigcup_{n \in \N} \bigcap_{m \geq n} a(m)) $, then there is $n \in \N$ and $y \in \bigcap_{m \geq n} a(m) $ such that $x \ll y$. As all $a(m)$ are past, $x \in a(m) $ for all $m \geq n$, thus $x \in \liminf_{2^X } (a)$. Let $x \in I^- (\limsup_{2^X} (a)) $, then there is $j: \N \nearrow \N $ with $x \ll a(j(n))$ for all $n \in \N$, and as all $a(j(n)) $ are past, $x \in a(j(n)) $ for all $n \in \N$, thus $x \in \limsup_{2^X} (a)$. The other items follow similarly. \hfill \qed

\begin{Theorem}
\label{IC}
\label{path-generated}
\begin{enumerate}
\item Let $Y$ be a chr. dense chr. space. If $I^\pm(p) $ is open $\forall p \in Y$, then $Y$ is regular.
\item  For a regular chr. space $(X, \ll, \tau)$ the maps $I_X^\pm: X \rightarrow IP(X)$ are inner-continuous.
\end{enumerate}

\end{Theorem}

{\bf Proof.} For the first assertion, we show that every $I^-(x)$ is synoptic (and thus indecomposable, see Sec.2): Let $y,z \in I^- (x) $, then $I^+(y) \cap I^+(z)$ is an open neighborhood of $x$, so chr. denseness shows that $I^-(x) \cap I^+(y) \cap I^+(z) \neq \emptyset $. For the second assertion, let $p \in X$ and $p \ll q$, then chr. separability ensures some $r_q \in I^+(p) \cap I^-(q) $, and regularity implies that $I^+(r_q) $ is an open neighborhood of $q$ whereas $I^-(r_q) $ is an open neighborhood of $ p$. Let $C \in I^+(p)$ be compact, then  the open covering $\{ I^+(r_q) \vert q \in C \} $ of $C$ has a finite subcovering $ \{ I^+ ( r_{q(i)}) \vert i \in \N_n \} $, and $U := \bigcap_{i=0}^n  I^-(r_{q(i)})$ is an open neighborhood of $p$ s.t. for each $q \in U$ we get $ C \in I^+(q) $. \hfill \qed

\bigskip

So, in a regular chr. space, $p \mapsto I^\pm(p) $ are inner continuous.  A future- or past-distinguishing chr. space $X$ is called {\bf future} resp. {\bf past causally continuous} iff $I^\pm$ is inner and outer continuous. A past subset of $X$ is determined by its closure: 
\begin{Lemma}
	\label{cl-inj}
 Let $A,B \subset X$ be past. If $\cl (A) = \cl (B) $ then $A= B$.
 \end{Lemma}
\V {\bf Proof.} Let $p \in A$, then $U:= I^+(p) \cap A \neq \emptyset$ is open; now $ U \subset \cl (A) = \cl(B) $, so $U \cap B \neq \emptyset$, thus $p \in B$ (and vice versa). \hfill \qed

\bigskip
\V For a chr. set $X$ we define the {\bf future} resp. {\bf past boundary} $\partial^+ X := \{ x \in X \vert I^+(x) = \emptyset\} $ resp. $\partial^- X := \{ x \in X \vert I^-(x) = \emptyset\} $. More generally, for a chr. space $X$ and $A \subset X$ we put $ \partial^+ (A, X):= \{ x \in \partial (A,X) \vert I_X^-(x) \cap A \neq \emptyset \}$. For conformal future-compact extensions, we have $ \partial^+E = \partial^+ (E(M), N) = \partial^+ X$ where either $X := \cl (E(M), N)$ or, restricting to the past-full part, $X := E(M) \cup \partial^+ (E(M))$, by the fact that both $\cl (E(M), N)$ and $E(M) \cup \partial^+(E(M), N) $ are chronological spaces chronologically embedded into $N$ (due to chr. convexity). In contrast to $E(M) \cup \partial^+ (E(M), N)$, the subset $\cl (E(M), N)$ is not semi-full due to the presence of $i_0$ as in the Penrose conformal compact extension of $\R^{1,n}$.

We say that a chronology $\ll$ on a topological space $X$ is {\bf path-generated} iff $p \ll q$ implies that there is a path (i.e., continuous map) $c$ from a real interval $I$ to $X$ with $c(0)= p$ and $ c(1)= q$ with $\forall s, t \in I: s <t \Rightarrow c(s) \ll c(t)$. Spacetime chronologies are path-generated. Of course, every path-generated chr. space is chronologically dense. By a little argument, each semi-full regular path-generated chr. space is locally path-connected (and recall that causal diamonds e.g. in Minkowski spacetime are not semi-full because of their spatial infinities).
If the chronology of a distinguishing chr. space $X$ is path-generated, $A \subset X$ is indecomposable if and only if $A$ is the past of the image of a chronological  future curve $c$ in $A$.

A regular path-generated chr. space $ (X, \ll , \tau) $ is called {\bf well-behaved} iff $(X, \tau)$ is locally compact, connected and locally arcwise connected (and recall that, by classical arguments using the pointwise compactness radius, each connected locally compact topological space is also sigma-compact).

Each chr. spacetime and each connected Lorentzian length space is well-behaved: From connectedness, local compactness and paracompactness we can conclude sigma-compactness, and chr. separability follows from localizability. Via Hausdorffness and chr. denseness we see that every well-behaved chr. space is $I$-distinguishing, i.e. for $I(A) := I^+(A) \cup I^-(A) $ we get that $I : X \rightarrow 2^X$ and $(I^+, I^-) : X \rightarrow 2^X \times 2^X$ are injective (whereas $I^\pm $ are not in general, due to possible nonempty causal boundaries). For the next theorem, let a chr. set $X$ be {\bf full} iff $I^+(x) \neq \emptyset \neq I^-(x)  \ \forall x \in X$. For a  chr. set $X$, we denote its largest full subset by $\underline{X}:= X \setminus (\partial^+ X \cup \partial^- X)$.

 \V A classical result by Vaughan (\cite{Vaughan}) says that each sigma-compact and locally compact topological space admits a compatible Heine-Borel metric. Even more, we can combine this with the result obtained in \cite{TT} that for each locally compact, connected, locally connected and separable metrizable space $(X, T)$ there is a $T$-compatible complete length metric on $X$ and incorporate it in Vaughan's proof to obtain a stronger statement:
 
 \begin{Theorem}
 	Let $X$ be a sigma-compact, locally compact, connected and locally connected metrizable topological space. Then there is a compatible intrinsic metric on $X$ that is Heine-Borel. 
 \end{Theorem}
 
\V{\bf Proof.} There is a sequence $U: \N\rightarrow 2^X$ of open relatively compact subsets with $\cl (U(n)) \subset U(n+1)$ for all $n \in \N$. Let $d$ be a compatible complete length metric on $X$ and

\( a_n := d(\cl (U(n)) , X \setminus U(n+1) ) = d(\cl (U(n)), \cl (U(n+2) \setminus U(n) )) >0 . \)   
 
 Here the equality holds as $d$ is a length metric and the inequality because of compactness of \\ $\cl (U(n+2) )\setminus U(n+1)) \subset \cl (U(n+2))$. For each $n \in \N$ we define $f_n:= \frac{1}{a_n} \cdot \max \{ d(\cl (U(n)), \cdot ), a_n \}$. Then $f_n$ is $a_n^{-1}$-Lipschitz continuous with $f_n \vert_{\cl (U_n)} = 0$, $f_n \vert_{X \setminus U(n+1)} =1$, and $g:= \sum_{k=1}^\infty f_k$ is well-defined and locally Lipschitz continuous, and the intrinsification $d''$ of $d': X \times X \rightarrow \R$ defined by $d'(x,y) := d(x,y) + \vert f(x) - f(y) | $ is finite and compatible by equivalence to $d \leq d''$: If $u: \N\rightarrow X$ converges to $v \in X$ w.r.t. $d$ then also w,r,t, $d''$ by local Lipschitzness of $g$. And $d''$ is Heine-Borel: Each $d''$-bounded subset is contained in some $U_n$ and thus has compact closure. \hfill \qed
 
 \bigskip

 For any well-behaved causally continuous chr. space we can pick such an intrinsic Heine-Borel metric $d_0$, which determines corresponding metrics $d_H$ on compact subsets and $d_1$ on closed subsets. Recall that for a length metric we always have $\cl (B(x.r)) = \ov{B} (x,r) := \{y \in X \vert d(x,y) \leq r \} \ \forall r >0$, and that the curve $c: [0; \infty) \rightarrow C(X)$ with $c(t) := \ov{B} (x,t)$ for each $t>0$ and $c(0) := \{ x \}$ is $d_1$-continuous. Theorem \ref{HB} tells us that the topology induced by $d_1$ (henceforth denoted $\tau_1$) does not depend on the metric chosen. The following theorem shows an interesting link between $d_1$-convergence and i.o.-convergence for almost past sets. We define $\ov{I}^\pm (p) := I^\pm (p) \cup \{ p \}$ for all $p \in X$. Let the set of open subsets of $X$ be equipped with the topology $\tau_{io}$ of inner and outer convergence. The map $\cap: (\tau(X), \tau_{io}) \times (\tau(X), \tau_{io}) \rightarrow (\tau(X), \tau_{io})$ is continuous. The corresponding statement for $d_1$-convergence is still true for almost past sets, cf. Item 2 of the following theorem:

\begin{Theorem}
\label{CSimpliesSC}
Let $X$ be a well-behaved chr. space. Let $d_0$ be an intrinsic Heine-Borel metric on $X$, and let $d_1:= d_1^{x_0, d_0}$ be the corresponding metric on the closed subsets as defined in Theorem \ref{Busemann}.
\begin{enumerate}
\item $\cl^* d_1$ is a metric on $\tau (X)$ generating $\tau_{io}$ on $P(X)$ for $\cl: P(X) \rightarrow CAP(X)$.	
\item $CAP(X)$ is closed in $(C(X),d_1)$.	

\item $\cap: (CAP (X) \times CAF(X), d_1 \times d_1) \rightarrow (C(X), d_1)$ is continuous.
\item $X$ is causally continuous $\Leftrightarrow I^\pm: (K(X) , d_1) \rightarrow (\tau_X, \tau_{io})$ is continuous 

$\Leftrightarrow \cl \o \ov{I}^\pm : (K(X), d_1) \rightarrow (C(X), d_1)$ is continuous $\Leftrightarrow \cl \o \ov{I}^\pm : X \rightarrow (C(X), d_1) $ is continuous\footnote{The necessity to use $\ov{I}$ instead of $I$ stems from the fact that we do not know fullness at the points and we have to avoid the empty set. In the causal setting explored in the last chapter this difficulty does not appear, as we can assume fullness for $J$.}.

\end{enumerate}
Furthermore, if $X$ is causally continuous, then it is {\bf locally chr. convex}, i.e., every neighborhood of a point $x \in X$ contains a chronologically convex subneighborhood of $x$.
\end{Theorem}

\V {\bf Proof.} To show that $\cl^* d_1$ is a metric it suffices to recall from Lemma \ref{cl-inj} that $\cl$ is injective on $P(X)$. For the first and second item, let $x \in A_\infty$. We want to show that $I^-(x) \subset \Int(A_\infty)$. As $x \in A_\infty$, there is $a: \N \rightarrow X$ with $a(n) \in A(n)  \ \forall n \in \N$ and $a(n) \rightarrow_{n \rightarrow \infty} x$. Let $z \ll x$, then by chr. separability there is $w \in I^+(z) \cap I^- (x)$, and $I^+(w)$ is an open neighborhood of $x$, so there is $N \in \N $ s.t. $\forall n \geq N : a(n) \in I^+(w)$. As $A_n$ is past, $w \in A_n  \ \forall n \geq N$, so $w \in A_\infty$, and $z \in I^-(w) \subset \Int (A_\infty)$.  

\medskip

For the second assertion of the first item, the statement on i.o.-convergence, let a compact $C \subset I^- (A_\infty)$ be given, then there is a finite covering of $C$ by subsets of the form $ I^-(z_i)$ with $z_i \in I^-(A_\infty)$, there are $y_i \in I^-(A_\infty) $ with $y_i \gg z_i$ for all $i \in \N_m$. Then for all $i \in \N_m$ there is a sequence $a^{(i)}$ with $a^{(i)} (n) \in A(n)$ for all $n \in \N$ and $a^{(i)} (n) \rightarrow_{n \rightarrow \infty} y_i$. By proceeding as before we show that there is $N \in \N$ such that $C \subset A(n) \ \forall n \geq \N$. Analogously for outer convergence, where we use that $\cl \o (X \setminus (\cl \o A))$ $d_1$-converges to $\cl (X \setminus \cl (A_\infty))$ (as the Hausdorff distance and thus also $d_1$ is compatible by construction with complements), and by Theorem \ref{AlmostPast}, $X \setminus A_n$ and $X \setminus A_\infty$ are almost future for all $n \in \N$ and $X \setminus I^-(A_\infty) = I^+(X \setminus A_\infty)$. 

\medskip

To prove the second item, for all $n \in \N$, let $A(n)$ be almost past, i.e., $I^-(A(n)) \subset A(n) $, and $B(n) $ be almost future, i.e., $I^+(B(n)) \subset B(n)$. We want to show that $d_1(A(n) \cap B(n), A_\infty \cap B_\infty) \rightarrow_{n \rightarrow \infty} 0$, or, equivalently by Th. \ref{HB}, for all $K \subset X$ compact and all $\e >0$ there is $N \in \N$ s.t. for all $n \geq N$: 

\bea
\label{Schnittgl1}
A_\infty \cap B_\infty \cap K &\subset&  B (A (n) \cap B (n) \cap K , \e)  \ {\rm and}\nonumber\\
 A(n) \cap B(n) \cap K &\subset&  B (A_\infty \cap B_\infty \cap K , \e).
\eea 

For the first inclusion, let $X= \{ x_1, ..., x_M\}$ be a finite $\e/2$-grid of $A_\infty \cap B_\infty \cap K $, i.e., $A_\infty \cap B_\infty \cap K \subset B(X, \e/2)$. For all $i \in \N_M^*$ we find points $y_i \ll x_i \ll z_i$ with $d(x_i, y_i), d(x_i, z_i) < \e/2$, then $Y:= \{ y_1, ... , y_M\}$ and $Z:= \{ z_1, ..., z_m \} $ are $\e$-grids of $A_\infty \cap B_\infty \cap K$. As $I^+(x_i)$ are open neighborhoods of $z_i \in A_\infty$ , we find $P_i \in \N $ s.t. for all $n \geq P_i$ we have $A(n) \cap I^+(x_i) \neq \emptyset$, and as $A(n)$ is almost past, we have $x_i \in A(n)  \ \forall n \geq P_i$. Analogously we find $Q_i \in \N $ suh that for all $n \geq Q_i$ we have $x_i \in B(n)$. Defining $N:= \max \{ \max \{ P_i \vert i \in \N_M^* \},  \max \{ Q_i \vert i \in \N_M^* \} \} $, the first inclusion of Eq. \ref{Schnittgl1} is true if $n \geq N$.  

The second inclusion in Eq. \ref{Schnittgl1} can be shown by quite general (non-causal) arguments: For all $\e >0$ we know, for large enough $n$, $A(n) \subset B(A_\infty, \e)  $ and $ B(n) \subset B(B_\infty, \e )$ and therefore $X \setminus A(n) \subset X \setminus B(A_\infty, \e) $ and $ X \setminus B(n) \subset X \setminus B(B_\infty, \e)$, thus the second inclusion follows from

\bean
X \setminus (A_n \cap \ B_n \cap K) &=& (X \setminus A_n ) \cup (X \setminus B_n) \cup (X \setminus K) \subset B(X \setminus A_\infty, \e) \cup B(X \setminus B_\infty, \e) \cup B(X \setminus K, \e) \\
&=& B((X \setminus A_\infty) \cup (X \setminus B_\infty) \cup (X \setminus K), \e) = B(X \setminus (A_\infty \cap B_\infty \cap K), \e),
\eean

For Item 3, we first number the statements whose equivalence is to be shown from (i) to (iv). We will show (i) $ \Rightarrow $ (ii) $ \Rightarrow $ (iii) $ \Rightarrow $ (iv), and we first observe that causal continuity is defined as continuity of $ I^\pm: X \rightarrow (\tau_X, \tau_{io}) $ and that (as the domain of every involved map is metrizable) it suffies to show sequential continuity in every case. 
(i) $ \Rightarrow $ (ii): Inner continuity: Let a sequence of compacta $K_n \rightarrow K_\infty$ in $d_1$ and let $L \in I^+(K_\infty)$. Then compactness of $L$ implies that there is a finite set $\{ p_1 , ... , p_n\} \subset K_\infty$ with $L \subset \bigcup_{k=1}^n I^+(p_i)$. By chr. denseness, for every $i \in \N_n^*$, there is $q^{(i)}: \N \rightarrow I^+(p_i)$ with $q^{(i)} (n) \rightarrow_{n \rightarrow \infty} p_i $, and causal continuity of $X$ implies that for each $i$ there is $n_i \in \N $ with $ L \subset \bigcup_{i=1}^n  I^+ (q(n_i)) $. As $I^- (q (n_i))$ is open for each $n_i$, there is $r_i$ such that $q(n_i ) \in I^+(K_M) $ for all $M \geq r_i$, then with $ r:= \max \{ r_1 , ... , r_n \} $ we get $q(n_i) \in I^+(K_M) \ \forall M \geq r_i  \ \forall i \in \N_n^*$, thus $L \subset I^+(K_M)$. For outer continuity proceed analogously: Let $S \in X \setminus I^+ (K_\infty) $ and let $ z_1, ... , z_n \in X $ with $U:=  \bigcup_{i=1}^n I^+(z_i) \supset K_\infty$ s.t. $S \cap I^+(z_i) = \emptyset$ for all $i$. Then there is $r \in \N$ s.t. for all $N \geq r$ we have $K_N \subset U$, and we conclude as above $L \cap I^+ (K_N) = \emptyset$.

(ii) $\Rightarrow$ (iii) follow directly from Theorems \ref{HausdorffInnerOuter} and Theorem \ref{HB}.

(iii) $\Rightarrow $ (iv) holds because for a sequence $a$ in $X$ convergent to some $a_\infty$, the sequence of compacta $ n \mapsto \{ a (n) \} $ converges in $d_1$ to $ \{ a_\infty \} $.   

(iv) $\Rightarrow $ (i): Let $a \in X^\N$ with $a(n) \rightarrow_{n \rightarrow \infty} p $. Item 1 then implies that $I^- (a(n)) \rightarrow_{n \rightarrow \infty}^{io} I^- (\cl (I^- (p))) =_{(*)} I^-(p),$
whose last equality (*) can be shown as follows: "$\supset$" holds as $ I^- (p) = I^-(I^-(p)) \subset I^- (\cl (I^-(p))) $, and "$ \subset$" holds as well: Let $x \in I^- (\cl (I^- (p)))$, i.e., there are $y \gg x$ and $z : \N\rightarrow I^-(p) $ with $z(n) \rightarrow_{n \rightarrow \infty} y$. As $I^+(x)$ is open, there is $n \in \N$ such that for all $N \geq n$ we have $z_N \in I^+(x)$, thus $x \ll p$.

For the first statement of the last item, let $U_1$ be a neighborhood of $x \in X$. Local compactness of $X$ ensures that there is a compact subneighborhood $U_2 \subset U_1$ of $x$. Metrizability and local compactness of $X$ imply that $X$ satisfies the $T_3$ separation axiom. Therefore there is an open neighborhood $U_3$ of $x$ with $V:= \cl(U_3) \subset \Int (U_2) $. The sets $V$ and $L:= U_2 \setminus U_3$ are compact. 

Now, defining $CB(x,r) := \{ y \in X \vert d(y,x) \leq r \}$ for $x \in X$, $r \geq 0$, let $ \Phi(r) := \cl (I^+(CB(x,r))) \cap \cl (I^-( CB(x,r))) $. For all $r >0$, $\Phi(r) $ is a neighborhood of $x$. As $\Phi (0) = \{ x \}$, we have $L \subset \Int (X \setminus \Phi(0))$, so there is $ r>0$ with $\Phi(r) \cap L = \emptyset$ due to outer continuity applied to $U_2$ (here we need that $r \mapsto (CB(x,r), d_1)$ is continuous at $r=0$ and the second and the third item). Local path-connectedness of $ X$ implies that $\Phi(r) \subset U_2$, thus $\Phi(r)$ is the desired neighborhood of $x$.

For the last assertion, the previous facts imply that there is $\e >0$ s.t. $U_\e := \cl (I^+ (B(x, \e)) ) \cap \cl (I^-(B(x, \e)))$ is contained in $U$ (as $U_0= \{ x \} \subset U$), it is chr. convex and contains $x$.  \hfill \qed

\newpage

\section{The topology $\tau_+$}
\label{MetricsOnIPs}

\v This section presents the topology $\tau_+$ on a chr. set $Y$ and applies it to the case $Y= IP(X)$ or $Y= P(X)$ for a chr. set $X$. Equivalence of the various definitions is shown in Th. \ref{TFAE}.

\bigskip

Let $(Y, \ll)$ be any distinguishing chr. set., let $a : \N\rightarrow Y$. We want to define convergence of $a$ in terms of $\liminf$ and $\limsup$ in an ordered set, and the obvious choice is $P(Y)$, to which we can map points in $Y$ by $i_Y = I^-_Y$. We define $a(n) \rightarrow_{n \rightarrow \infty} a_\infty$ iff 

\[ \limsup \ _{P(Y) } I^-_Y  \o a = \liminf \ _{P(Y)} I^-_Y \o a = I^-_Y (a_\infty). \]

 Clearly $\limsup (a) \geq \liminf (a)$, as the union in the definition of the latter is not over all subsequences but only over those of the form $ n \mapsto c+n$ for some constant $c \in \N$. Furthermore, for a monotonically increasing sequence $\limsup$ and $\liminf$ coincide. We define $L_+: Y^\N \rightarrow 2^Y $ by
 
 $$ L_+(a) := \{ v \in Y \vert I^-(\liminf \ ^- (a)) = I^-(\limsup \ ^- (a)) = I^-(v) \} $$
 
for every $a \in Y^\N$, and $S \subset Y$ is $\tau_+$-closed iff for every sequence $a$ in $S$ we have $L_+(a) \subset S$. $L_+$ is called {\bf limit operator of $\tau_+$}. As $Y$ is distinguishing, we get $ \# (L_+(a)) \in \{ 0, 1 \} $ for all $a: \N\rightarrow Y$.

This indeed defines a topology $\tau_+$ on $Y$, due to a classical construction by Fr\'echet and Urysohn (see e.g.  \cite{Engelking}, p. 63): We only have to show that for all $a: \N \rightarrow Y$:

\begin{enumerate}
\item If $a(\N ) = \{ p \}$ we have $L_+(a) = \{ p \}$;
\item If $L_+(a) = \{ p \}$ then for every subsequence $b$ of $a$ we have $L_+(b) = \{ p \}$;
\item If $p \notin L_+(a)$ then there is a subsequence $b$ of $a$ s.t. every subsequence $c$ of $b$ has $ p \notin L_+(c) $.
\end{enumerate} 

These conditions are easy to verify: The second one is a consequence of $ \liminf^- (a) \leq \liminf^-  (b) \leq \limsup^- (b) \leq \limsup^- (a) $ for any subsequence $b$ of a sequence $a$ (recall that $\leq = \subset$ on $P(Y)$). To prove the third one, either we have $I^-(\limsup^- (a)) \neq I^-(\liminf^- (a) )$ in which case we choose $b$ as a subsequence making up this difference, i.e. with $\liminf^- a \neq \bigcup_{n \in \N } I^{\pm \bigcap} (b(\N)) $; or $I^-(\limsup^- (a)) = I^-(\liminf^- (a)) $ but then this holds for all subsequences, which finishes the argument. So indeed $\tau_+: C \rightarrow CS$ is a well-defined functor, a left inverse of the forgetful functor.

\medskip

There are strongly causal spacetimes whose convergence structure is not described correctly by $L_+$: Let $\R^{1,1} := (\R^2 , - dx_0^2 + dx_1^2)$ and consider $M:= \R^{1,1} \setminus ( \{ 0 \} \times [0; \infty))$, a strongly causal spacetime that is not past-reflecting. For $z(n) := (1+1/n, 1) $ and $z_\infty := (1,1)$, the sequence $z$ converges to $z_\infty$ in the manifold topology but not in $L_+$, as $ A:= I^-_{\R^{1,1}} (0) \subset I_M^-(z_n)  \forall n \in \N$ but $ A \cap I_M^-(z_\infty) = \emptyset$. Theorem \ref{recovering} will show that this discrepancy only happens in the non-causally continuous case.

\medskip

We want to connect the previous notions of limsup and liminf to the set-theoretical ones in the case of a chr. set $Y= IP(X)$, where $X$ is a regular and $I^-$-distinguishing chr. space. To this aim, we consider the relation $\ll_{BS}$ on $2^X$ and the sets $I^\pm_{2^X} $ defined by it\footnote{Here we just infer the definition $I^-(A) := \{ B \in 2^X \vert B \ll_{BS} A \}$ but we recall that $\ll_{BS}$ is not a chr. relation on $2^X$ in general, as chr. separability might fail.}. This is necessary, as $I^-_X (\liminf (a))$ is not indecomposable in general for a sequence $a$ of indecomposable subsets: Consider, e.g., the sequence $ a: \N \rightarrow IP(\R \times \mathbb{S}^1) $ into the two-dimensional Einstein cylinder given by $a(n) := (0,(-1)^n )$ for all $n \in \N$, then $\liminf (a) $ is the union of two PIPs.

\begin{Theorem}[see also \cite{CFH}]
Let $X$ be a past-full chr. set. On $IP(X)$, we have $\tau_- \subsetneq \tau_+$.
\end{Theorem}

\V {\bf Proof.} Let $ a : \N \rightarrow IP(X)$. The synopsis of the Definitions 2.3, 2.4 and 6.1 from \cite{jF} shows that $L_-(a) $ is precisely the set of those $v \in X$ such that for every $A \subset X$ that is a maximal IP (resp. maximal IF) in $I^-(v)$ (resp. $I^+(x)$) satisfy: $A \subset \liminf a $ and $ A$ is a maximal IP (resp. IF) in $\limsup a $. The definition of $ L_+   $ on $IP(X) $ shows that $L_+ \subset L_-$ as partial maps from $(IP(X))^\N$ to $P(IP(X))$ (because $I^- (\liminf^- (a) ) = I^-(v)$ implies $A \subset I^-(x ) \Rightarrow A \subset \liminf (a) $ and $I^-(\limsup^- (a)) = I^-(v) $ implies that if $A$ is a maximal IP in $I^-(x) $, then $A$ is a maximal IP in $\limsup (a)$), which implies the claim. The fact that the two topologies do not coincide in general, even on globally hyperbolic spacetimes, is shown by Harris' example of the unwrapped-grapefruit-on-a-stick \cite{H} accounted for at the end of this section. \hfill \qed

\bigskip

The next task is to examine $\tau_+$ on $IP(X)$. First, the definition of the metrizable topology in \cite{oM14} was made for $IP(M)$ where $M$ is a g.h. spacetime. In this case, when defining an appropriate measure, we can use Lemma 3.3 of \cite{CGM} stating that, for a globally hyperbolic manifold $(N,h)$ and for any compactly supported $\psi \in C^0 (N, [0, \infty))$ with $\int_M \psi (x) d {\rm vol} (x) = 1$, the function $\tau_{\psi}$ with $\tau_\psi (p) := \int_{I^-(p)} \psi d{\rm vol}_h$ is continuously differentiable. We choose a locally finite countable covering of $M$ by open precompact sets $U_i$ and define, for $\phi_i \in C^ \infty (M, (0, \infty))$ with $ \phi_i^{-1} (0) = M \setminus U_i $, 

$$ \phi = \sum_{i \in \N} 2^{-i} (\vert \vert \phi_i \vert \vert_{C^1} + \vert \vert \phi_i \vert \vert_{L^1})^{-1} \phi_i , $$ 

therefore $\phi \in C^1 (M) \cap L^1 (M) $ with $\vert \vert \phi \vert \vert_{C^1(M)} , \vert \vert \phi \vert \vert_{L^1(M)} \leq 1$; we rescale $\phi$ such that $\vert \vert \phi \vert \vert_{L^1(M)} =1 $.

When treating general chronological spaces $X$, we have to ensure the existence of appropriate measures on them: A synopsis of \cite{Knowles} (Corollary after Th. 3) and \cite{HebertLacey}(Cor.2.8) ensures that on every Polish (i.e., complete-metrizable separable topological) space without isolated points there is an {\bf admissible} measure, i.e. a finite non-atomic strictly positive Borel measure (the finiteness it not mentioned in the corollary but clear from the construction via the Stone-Cech compactification). A classical result \cite{Bauer}(Lemma 26.2) ensures that every finite Borel measure $\mu$ on a Polish space is Radon, that is, locally finite, outer-regular and inner-regular: for every Borel set $B$ we have $\mu (B) = \sup \{ \mu (K) \vert K \subset B {\rm \ compact} \} = \inf \{ \mu (U) \vert B \subset U {\rm \ open} \}$.
In summary, on every Polish space $X$ we have a non-atomic strictly positive  finite Radon Borel measure $\mu$, which we can apply instead of $\Phi \vol$ above. Now, we can induce a metric $\delta_\mu : IP(X) \times IP(X) \rightarrow [0; \infty) $ by defining

$$\delta_\mu (A,B) :=  \mu (\Delta (A,B)).  $$

\V where $\Delta(A,B) := (A \setminus B) \cup (B \setminus A)$ is the symmetric difference of $A$ and $B$.

For the last statement of the following theorem, assuming that $X=M$ is a g.h. spacetime, we use a Cauchy temporal function $t$ (whose existence is ensured by the nowadays classical result \cite{BS}, for a short self-contained account and extended results see \cite{oM15}) to define a diffeomorphism $D: \R \times S \rightarrow M$ with $ D^* t = {\rm pr}_1 $ (where $S$ is a Cauchy surface of $M$) and $D^* g = - L^2 dt^2 +  {\rm pr}_2^* (h \o t) $ where $L$ is a function on $\R \times S$ and $ h: \R \rightarrow {\rm Riem} (M) $ is a one-parameter family of Riemannian metrics on $S$. As for every IP $A$ in $(\R \times S, D^* g)$, $ \partial A$ is an achronal boundary, which, as a subset of $\R \times S$, is well-known (see e.g. \cite{eM}, Th. 2.87(iii)) to be (the graph of) a locally Lipschitz partial function $f (t,A): S \rightarrow \R$. The Hausdorff metric w.r.t. a metric $g$ is denoted by $d_{H(g)}$. For a set $D$, let ${\rm Metr} (D)$ denote the set of metrics on $D$. Recall from Th. \ref{CSimpliesSC} that the pullback $\cl^ *d_1$ of $d_1 $ along $\cl$ is a metric on $P(X)$.

\begin{Theorem}
\label{TFAE}
Let $X$ be a well-behaved past-full chr. space, let $a \in (P (X))^\N $, let $a_\infty \in P(X)$. Let $d$ be a compatible Heine-Borel metric on $X$. Then the following are equivalent:

\begin{enumerate}
\item $\cl \o a$ converges in $(C(X), d_1)$ to $\cl (a_\infty)$,
\item For every finite non-atomic strictly positive Borel measure $\mu$ on $X$ and for every compact subset $K$ of $X$, $\lim_{n \rightarrow \infty} (\mu (\Delta (a_\infty \cap K, a(n) \cap K) )) = 0$ and $\d_\mu: IP(X) \times IP(X) \rightarrow \R$ defined by $\delta_\mu (A,B) := \mu (\Delta (A,B)) $ is a metric topologically equivalent to $d_1$,
\item $\forall C \in K(X) :$ $\lim_{n \rightarrow \infty} (d_{H (g)} (\cl (a_\infty \cap C), \cl (a(n) \cap C)))  = 0 \ \forall g \in {\rm Metr} (C)$ compatible,
\item $I^- (\limsup a) = I^-(\liminf a) = a_\infty$,
\item $\Int (\limsup (a)) = \Int (\liminf (a)) = a_\infty  $,
\item $\cl (\limsup (a)) = \cl (\liminf (a)) = \cl (a_\infty)$,
\item $\widetilde {\limsup (a)} = \widetilde {\liminf (a)} = \cl (a_\infty)$,
\item $L_+(a) = \{ a_\infty\}$,
\item $a$ converges in $(P(X), \tau_+)$ to $a_\infty$.
\end{enumerate}

In particular, 8. $\Leftrightarrow$ 9. shows that $(P(X), \tau_+, L_+) $ is of first order.
If $X$ is a g.h. spacetime, all nine conditions are equivalent to: For $t$ a continuous Cauchy time function and $f(t,A)$ defined as above, directly before this theorem, $f (t,a(n))$ converges to $f(t,a_\infty)$ pointwise, or equivalently, in the compact-open topology (*).
 \end{Theorem}
 
\V {\bf Remark.} Beem's original definition of $\tau_+$ in \cite{Beem!} is Item 7. In his Proposition 3 he shows equivalence of Items 1 and 7, and he states without proof 1. $\Leftrightarrow$ 2. for manifolds with compact Cauchy surface.

\V {\bf Remark.} \v Unlike the topology $\tau_+$, neither the metrics inducing it nor their uniform structures are functorial, by the unfunctorial choice of $\mu$. But we can construct natural metrics e.g. in the category of temporally compact g.h. spacetimes, cf. \cite{oM-Finiteness}.

 \medskip

\V {\bf Proof.} First we want to show that Items 1,3,4,6,7 and 8 are equivalent without further conditions.

$(1) \Leftrightarrow (3)$: see Th. \ref{HB}.
%$ d_1 (a_\infty , a(n)) > d_H ( (a_\infty \cap K , a(n) \cap K))$. 

$(1) \Rightarrow (7)$: This is the second last assertion of Th. \ref{Busemann}.

$ (1) \Rightarrow  (5)$: Let $p \in a_\infty$, then the open subset $A:= I^+(p) \cap a_\infty$ of $X$ contains a nonempty open subset $U$ with $ K= \cl(U,X) \subset A $ compact. By the assumption (1), there is $n \in \N$ such that for all $N \geq n$ we have $a_N \cap U \neq \emptyset$, implying $p \in a_N$. Consequently, $p \in \liminf (a)$. Conversely, if there is $n \in \N$ such that for all $N \geq n$ we have $p \in a_N$, then $ p \in \cl (a_\infty)$: Assume $p \in X \setminus \cl (a_\infty)$, then for an open and precompact neighborhood $B$ of $p$ whose closure is contained in $ X \setminus \cl (a_\infty) $, we have $B:= I^-(p) \cap K \subset X \setminus \cl (a_\infty)$; on the other hand, $B \subset a_n  \forall n \in \N$ as the $a_n$ are past; thus if a ball of radius $r$ is contained in $U \subset B_R(x_0)$ then $ d_1 (a_\infty , a_n) \geq  r e^{-R} >0 \ \forall n \in \N$, in contradiction to (2). Thus $\liminf (a) \subset \cl (a_\infty)$, so $\Int (\liminf (a) ) \subset a_\infty$, and with the above $a_\infty = \Int(\liminf (a)) \subset \Int (\limsup (a)) $.

To show $\Int (\limsup (a)) \subset a_\infty$: First, $a_\infty \in P(X)$ is open, secondly $\Int (\limsup (a) ) \ \subset \cl (a_\infty)$:

Let $p \in \Int (\limsup (a))$, then there is a strictly increasing $j: \N \rightarrow \N$ with $p \in \bigcap_{k=0}^\infty a(j(k))$. We want to show $p \in \cl (a_\infty) $ by constructing a sequence $y: \N \rightarrow a_\infty$ with $\lim (y) = p$. As $X$ is past-full, there is a sequence $b$ in $I^-(p) \cap \Int (\limsup a) $ converging w.r.t. $d_1$ to $p$. For all $n \in \N$ we define $U_n:= I^+(b(n)) \cap I^-(p)$. We have a ball of radius $r_n>0$ contained in  $U_n \in B(x_0, R)$ for all $n$, so 

$$  \forall n \in \N  \exists l \in \N \forall L \geq l : d_1 (a(L), a_\infty ) < r_n e^{-R}. $$   

As, on the other hand, $U_n \subset a(j(k)) $ for $j(k) > L$, we have $\emptyset \neq a_\infty \cap U_n \ni y_n $, defining a sequence $y$ with $y(n) \rightarrow_{n \rightarrow \infty}p$ via the "causal sandwich principle" in $(P(X), d_1)$ implied by Th. \ref{Busemann}: If $p \in (I^-(p_\infty))^\N$ converges to $p_\infty$ and $q \in X^\N$ with $p(n) \ll q(n) \ll p_\infty \forall n \in \N$ then $q(n) \rightarrow_{n \rightarrow \infty} p_\infty$.

$(4) \Leftrightarrow (8)$: by definition.

$(6) \Leftrightarrow (7)$: For sequences $a$ with $J^-  \o a = a$ we have $ \cl (\limsup (a)) = \widetilde{\limsup} (a) $ and $  \cl (\liminf (a)) = \widetilde{\liminf} (a )$: Let $p \in \widetilde{ \liminf  } (a) $. Chronological denseness implies that there is $ k: \N \rightarrow I^-(p)$ (w.l.o.g. chosen as a chr. chain) converging to  $p $. Openness of the $I^+(k(n))$ and the condition that the $a(n) $ are past imply that there is a subsequence $l$ of $k$ with $ l(n) \in a(n) $. The reverse inclusion is trivial.

$ (7) \Leftrightarrow (1) $ follows immediately from the fact that $d_1$-convergence is equivalent to the topology determined by convergence of $\widetilde{\limsup}$ and $\widetilde{\liminf}$ (second last assertion in Th. \ref{Busemann}) and the fact that $ \widetilde{\limsup} ( \cl \o a ) = \widetilde{\limsup} (a)$ (and corespondingly for $\widetilde{\liminf})$ as $\cl \o \cl = \cl$. 

$ (4) \Rightarrow (6)$: By Th. \ref{AlmostPast}, $\limsup (a) \in AP(X)$, so $\cl (\limsup (a) ) = \cl (I^- (\limsup a))$, also for $\liminf$.

$(5) \Leftrightarrow (6)$: We apply $\Int$ resp. $\cl$ to each side and use Theorem \ref{AlmostPast}.

$(5) \Rightarrow (4)$: Apply $I^-$ to (5), use $\limsup a) \in AP(X)$, $I^-(\Int A) = I^-(A)$ for $A \in AP(X)$ (Th. \ref{AlmostPast}). 

Thus we get $(1) \Rightarrow (4) \Rightarrow (6) \Leftrightarrow (7) \Leftrightarrow (1) \Leftrightarrow (3) $, showing pairwise equivalence of (1),(3),(4),(6),(7), (8). As the topology in (1) is metric, the limit operator is first order and we get equivalence to (9).  

$(2) \Leftrightarrow (3)$:  First of all, any $\sigma$-compact topological space is separable. Also any $\sigma$-compact metrizable topological space is complete-metrizable. Consequently, every well-behaved space is Polish. Now, for the first assertion, $\delta_\mu$ is obviously nonnegative and symmetric. The triangle inequality follows immediately from the set-theoretic triangle inequality $\Delta(A,C) \subset \Delta(A,B) \cup \Delta (B,C) $. It remains to show that $\d_\mu$ does not vanish between different subsets. By regularity, past subsets are always open. But in general, for $A$,$B$ two different open past subsets, we have $(A \setminus \overline{B}) \cup (B \setminus \overline{A}) $ is nonempty and open: Let, w.r.o.g., be $x \in A \setminus B$. Openness of $A$ allows to choose a nonempty open precompact $U \subset X$ with $ \cl( U ) \subset A \cap I^+ (x)  $.  By the fact that $B$ is past and does not contain $x$ we know that $U \cap  \overline{B} = \emptyset$. Consequently, $U \subset A \setminus B $, and $ \mu (\Delta(A,B))   \geq \mu (U)  >0 $. Thus $\delta_\mu$ is not only a pseudometric but a true one. Now, the claim follows from inner and outer regularity of $\mu$ which implies $d_1$-continuity of $\mu$ on the set of closures of past subsets according to Theorem \ref{CSimpliesSC}.

$(9)$ is equivalent to $(8)$ due to the equivalence of $(8)$ and $(1)$: The limit operator $L_+$ is equal to the limit operator w.r.t. the metric $d_1$, therefore the closed sets are the same, thus the topology.

$(4) \Leftrightarrow (*)$ for g.h. spacetimes (cf. \cite{CFH}, Th. 3.9) follows from the fact that for locally Lipschitz maps pointwise convergence is equivalent to compact-open convergence on a compact set. Here it is more practical to show $(*) \Leftrightarrow (3)$: Assume (3), let $C \subset X$ be compact, let $K:= C \cap t^{-1} (\{ 0 \})$. Choosing a compatible a product metric, with $f_n:= f(t, A_n)$ and $f_\infty:= f(t,A_\infty)$, we get $f_\infty \vert_K$ and $f_n \vert_K $ Lipschitz continuous for all $n \in \N$, with the same Lipschitz constant. Due to compactness of $K$ and Lipschitz continuity, for every $\e >0$ there is $\d>0$ s.t. for each $p \in K$ with $\vert f_n(p) - f_\infty(p) \vert > \e$ there is a ball around $(p, f_n(p)) \in \partial a_n $ of radius $\d$ not intersecting $a_\infty$ (if $f_n(p) > f_\infty (p)$) or there is  a ball around $(p, f_\infty(p)) \in \partial a_\infty$ of radius $\d$ not intersecting $a_n$ (if $f_n(p) < f_\infty (p)$). In either case, $d_H (C \cap \cl(a_n), C \cap \cl(a_\infty) ) \geq \d$. So $(3) \Rightarrow (*)$. The other implication is the well-known fact that for compact $K$, if $f_n \rightarrow f_\infty$ uniformly, $d_H (f_n \vert_K, f_\infty \vert_K) \rightarrow_{n \rightarrow \infty} 0 $ as subsets of $K \times \R$. \hfill \qed

\bigskip

The previous theorem shows that the topologies from \cite{oM14} and \cite{CFH} coincide. The concrete definition in \cite{oM14} is appropriate e.g. for calculations of future boundaries of spacetimes as those done in \cite{oM14}, whereas the original definition of $\tau_+$ due to Beem shows that the topology is functorial in the purely chr. category.

\begin{Theorem}
\label{recovering}
Let $(X,  \ll, \tau )$ be a causally continuous well-behaved past-full chr. space. Then $\tau = i_X^{-1} (\tau_+)$, in particular $i_X(X)$ is open and $i_X: X \rightarrow (i_X(X), \tau_+)$ is a homeomorphism. Moreover, $(P(X), \tau_+)$ and $IP^+ (X) := \cl (i_X(X), (P(X), \delta))$ (with induced chr. structure and topology) are again past-full causally continuous well-behaved chr. spaces both containing $IP(X)$. The relative topology of the subset $\underline{IP(X)}$ of full points of $(IP(X), \tau_+)$ is the Alexandrov topology. Finally, $i_X (c(n)) \rightarrow_{n \rightarrow \infty} I^- (c(\N)) $ for any chr.~chain $c$ in $X$.  
\end{Theorem}  

{\bf Remark:} So $\tau_+$ applied to a causally continuous spacetime recovers the manifold topology. Moreover, Theorem \ref{recovering} allows to define a functor $F_+ = (IP^+, \tau_+)$ from the category  $CW$ of causally continuous well-behaved past-full chr. spaces to itself, and $F_+$ is marginal by definition.

Furthermore, the theorem implies that $\tau_+$ is second-countable on $P(X), IP^+(X)$ and $IP(X)$ for $X$ well-behaved (for $X$ spacetime this is in Th. 5.2. (ii) of \cite{FHS1}), as metrizable and sigma-compact topological spaces are second countable (compact metric spaces are second countable, and a countable union of countable sets is countable). 

\bigskip

\V {\bf Proof.} First we want to show {\em regularity} of $i_X^{-1} (\tau_+)$, i.e. that $i_X (I^+(x) ) $ is $\tau_+$-open for all $x$. Given metrizability and thus sequentiality of $\tau_+$, it is enough to show that $I^+(x)$ is sequentially open, or, equivalently, that $X \setminus I^+(x)$ is sequentially closed. Therefore, let $a \in (X \setminus I^+(p))^\N $ converge to some $v \in X$, and we have to show that $ v \in X \setminus I^+(p)$. As $a(n) \in X \setminus I^+(p)  \forall n \in \N$, we know $p \not\ll a(n) \forall n \in \N$. Convergence of $a$ to $v$ means $I^\pm (\liminf^\pm (a)) = I^\pm (\limsup ^\pm (a)) = I^\pm (v)$, and we want to show $p \notin I^-(v)$, which with the convergence condition is equivalent to 
$p \notin I^-(\liminf^- a ) = I^- ( \bigcup_{n \in \N} I^{- \cap } (\{ a(m) \vert m \geq n \})) $. Assume $p \in I^- ( \bigcup_{n \in \N} I^{- \cap } (\{ a(m) \vert m \geq n \}))  $, then there is $n \in \N$ such that for all $m \geq n$ we have $p \in I^- (I^- (a(m))) = I^- (a(m))$, which is a contradiction.

 {\em $IP(X) $ is path-generated}: Let $P \ll Q \in IP(X)$, then either $P,Q \in i_X(X)$, in which case $P:= i_X (p), Q:= i_X(q)$ for some $p,q \in X$ with $p \ll q$. Then there is a chronological path $c: p \leadsto q $, and $i_X \o c$ is a chronological path from $P$ to $Q$. If, on the other hand, $P = i_X (p), Q \in \partial^+ X$, then $Q:= I^-(c) $ for a chronological path $c: (a;b) \rightarrow X$, i.e. there is $s \in \R$ with $p \in I^-(c(s))$, so there is a chronological path $k_1 : p \leadsto c(s)$ and for $k_2:= c \vert_{[s, b)}$, $ k:= k_2 * k_1$, $C: t \mapsto I^- (k (t))$ is a chronological path from $P$ to $Q$.

By Th. \ref{Busemann}, $(C(X),d_1)$ is locally compact (if $(X,d)$ was locally compact). Theorem \ref{TFAE}, equivalence of Items 1 and 9, implies that the topology $\tau_+= \cl^* (\tau_{d_1})$ on $P(X)$, generated by $D = \cl^* d_1$, also is locally compact, as $\cl $ is a homeomorphism.

{\em Metrizability} has been shown in Th. \ref{TFAE}. If $X$ is locally arcwise connected, w.l.o.g. consider a precompact neighborhood $U$ at $\partial^+ X $ contained in some $I^+(q)$.

The {\em assertion on $\underline{IP} (X)$} follows directly from the last assertion (local chr. convexity) in Th. \ref{CSimpliesSC} and chr. denseness.

\begin{Lemma}[Non-imprisonment]
	\label{LCC}
In a well-behaved locally chr. convex chr. space $X$, every chr. chain $c$ in a compact subset has an endpoint $p_c$, i.e. with $I^-(p_c) = I^-(c) $.
\end{Lemma}

\V {\bf Proof of the lemma}. The causal sandwich principle of Th. \ref{CSimpliesSC} implies that in the full subset $\underline{X}$ there is a chr. convex neighborhood basis (i.e., $\underline{X} $ is locally causally convex), by which one shows easily that an accumulation point is actually a limit (as each chr. chain cannot leave and re-enter a chronologically convex set).

\bigskip

\begin{Lemma}[Continuity of $i_X$]
Any sequence $x: \N \rightarrow X$ that $\tau$-converges to $x_\infty \in X$ also converges to $x_\infty$ in $i_X^*(\tau_+)$.
\end{Lemma}

{\bf Proof of the lemma.} Directly from causal continuity and Th. \ref{CSimpliesSC}, Item 4. 

\bigskip

\begin{Lemma}[Openness of $i_X$]
Let $p : \N \rightarrow (IP(X), \tau_+)$ converge to a point $I^-(x) \in i_X(\underline{X})$. Let $U$ be a $\tau$-open precompact neighborhood of $x$. Then there is $m \in \N$ with $ \forall n \geq m: p(n) \in i_X(U)$.  
\end{Lemma}

{\bf Proof of the lemma.} As seen in Lemma \ref{LCC}, $\underline{X}$ is Alexandrov, thus we find $p_{\pm} \in X$ s.t. $A:= \cl (I^+ (p_-) \cap I^-(p_+) ) \subset U$ is an compact neighborhood of $x$, and two other points $p'_\pm$ with 

\bea
\label{causal-ordering}
p_-  \ll p'_- \ll  x \ll p'_+ \ll  p_+ .
\eea

Eq. \ref{causal-ordering} implies the existence of $N_1, N_2 \in \N$ with 

\bea
\label{Folgerungen1}
\label{Folgerungen2}
p'_- \in \lim_{n \rightarrow \infty} p_n \subset \bigcap_{n \geq N_1} p_n, \\
p'_+  \in  \bigcap_{n \geq N_2} \cl (I^{+ \cap} (p_n)).
\eea

(last assertion due to openness of $I_{IP(X)}^+(x)$). Thus there is $N:= \max \{ N_1, N_2 \} \in \N$ with $p_- \in \cl (I^-(p_m)) \land p_+ \in \cl (I^{+ \cap} (p_n)) \ \forall n \geq N$, so for all $n \geq N$ we have $p_n= I^-(c_n)$ for a chronological chain $c_n : \N \rightarrow \cl (I^+ (p_-) \cap I^- (p_+)) = A $. As $ A $ is a compact subneighborhood of $U$, and $ X$ is nonimprisoning as shown in the Lemma above, the $c_n$ converge to some point in $U$. Thus, as we desired to show, $p_n = I^- (x_n)$ with $x_n  \in U$.

This showsthe first assertion of the theorem. {\em Local compactness} of  $(IP(X), \tau_+)$ follows as $(C(X), d_1)$ is locally compact and $\cl: (P(X), D= \cl^* d_1) \rightarrow (CAP(X), d_1)$ is a homeomorphism.

As for each $p \in \partial^+X $ there is a chr. chain approaching it and $X$ is arc-connected, the chain can be made a curve and $X$ is arc-connected.
Finally, to show {\em local arc-connectedness} of $IP(X)$, let $U$ be an neighborhood of $p \in \partial^+X$. As $X$ is locally compact, there is a precompact open neighborhood $V \subset W$ of $p$. As $\partial^+X$ is the set of limits of chr. chains, there is $q \in X$ with $I(q,p) \subset V$. We have $\cl (I(q,p)) = \cl (I^+(q) ) \setminus I^+(C)$ for the compact set $C:= \partial^+ (I^-(p) \cap I^+(q)) \subset V$. As $C$ is compact and contained in the open set $V$, $\rho:= d_1 (C, X \setminus V) >0$. Define $C^+:= I^+(C) \cap B(C, \rho/2)$. Then $W:= I^+(q) \setminus \cl (I^+(C^+))$ is an open neighborhood of $p$ contained in $V$ and thus in $U$. It is arc-connected: For each two $x,y \in W$, there are past curves $x \leadsto q$ and $y \leadsto q$.  \hfill \qed

\bigskip

If a topology $\tau $ on $X$ is causally continuous, then obviously every $\tau$-converging net $\tau_+$-converges, thus $\tau_+$ is the coarsest causally continuous topology. 

\bigskip
\V The next theorem (in parts adapted from \cite{oM14}, compare also \cite{FHS1}, Th. 4.16 for the corresponding assertion on $\tau_-$) shows that the intrinsic future completion defined above is homeomorphic to the union of the image of a conformal extension $E$ of a g.h. spacetime $M$ and its future boundary. The map $\e_E$ is defined as follows: Let $A \in IP(M)$, then choose a chronological chain $c: \N \rightarrow M$ with $A = I^-(c(\N))$. We define $\e_E(A)$ to be the limit of $E \o c$ in $\cl(E(M), N)$ (which exists due to future compactness of $E(M)$ in $N$). As argued before, this does not depend on the choice of $c$. Moreover, for the statement below it makes no difference which Cauchy surface $S$ one chooses in the definition of a CFE or if we even replace $I^+ (S)$ with $M$. For the definition of $\underline{X}$ see the second last paragraph before Th. \ref{CSimpliesSC}.

\bigskip

The relation between $\tau_-$-convergence and $\tau_+$-convergence is completely described in the following theorem:

\begin{Theorem}
	\label{ConvergencePlusMinus}
	 $(P(X), \cl^* d_1)$ is future-compact, i.e. $I^+(C)$ is precompact for all compact $C \subset P(X)$. For $p \in X$, each sequence $a: \N\rightarrow I^+ (p) $ has a subsequence $b= a \o j$ convergent in $(P(X), \cl^*d_1)$ to some $A \in P(X)$ which is, in general, decomposable. Let $U$ be a maximal indecomposable past subset in $A$. Then $b (n) \rightarrow^{\tau_-}_{n \rightarrow \infty} U$.
	 \end{Theorem}

\V{\bf Proof.} Future compactness follows from the fact (from the last assertion in Th. \ref{Busemann}) that for each $p \in X$, the set $ C_p := \{ A \in C(X) \vert p \in A \}$ is compact in the metric $d_1$. By Theorem \ref{TFAE}, $\cl (b(n)) \rightarrow_{n \rightarrow \infty}^{d_1} A $ implies $ A = \cl (\limsup (a(n))) = \cl (\liminf (a(n)))$, which, with the definition of $\tau_-$,  implies $b (n) \rightarrow^{\tau_-}_{n \rightarrow \infty} U $ for every maximal indecomposable past subset $U$ of $A$.  \hfill \qed

\bigskip

The previous theorem shows that $(P(X), \tau_+)$ contains at least as much information about the chr. structure of $X$ as $(IP(X), \tau_-)$. Actually, the former topology contains {\em strictly more} information than the latter, which can be seen in the example of $X := I^- (c((0; \infty)))$ for $c: (0; \infty) \rightarrow \R^{1,2}$ given by $c(t) := (t, \sin (1/t), 0)$: Consider the sequence $a: \N\rightarrow \partial^+X $ given by $c(n) := I^- (c(1/n))$, then $a$ does not converge in $\tau^+$, as $\limsup (a) = I^- (\{ 0\} \times [-1;1] \times \{0 \}) $ but $\liminf (a) = x_3^{-1} ((- \infty; - 1))$. In $\tau_-$, however, the sequence converges to each $I^- ((0,t,0))$ with $t \in [0;1]$, whereas, in contrast, we will see an example (due to Harris) of a sequence displaying the same behaviour (convergence to more than one element of $IP(X)$) w.r.t. $\tau_-$ while not converging at all in $P(X) $ with $\tau_+$.

\bigskip

The previous theorem also entails an interesting aspect of the causal boundary: it admits to determine maximal and minimal elements (w.r.t. the order given by inclusion) in the future of a compact subset\footnote{recall that an element $x$ of an ordered set $(X, \leq)$ is called {\bf maximal} iff $ \forall y \in X : \neg (x < y)  $, correspondingly for minimal elements. On $IP(X)$ inclusion corresponds precisely to causal relatedness. More on this in Sec. \ref{posets}.}.

\begin{Theorem}
	Let $C \subset IP(X)$ be compact, then there are minimal and maximal elements in $J^+(C) \cap \partial^+X$.
\end{Theorem}

{\bf Remark.} This problem has brought to the author's attention by Leonardo Garc\'ia-Heveling. 

\bigskip

\V{\bf Proof.} Consider a $\mu$-infimizing sequence $a$ in $J^+(C)$, i.e., $\mu (a_n) \rightarrow_{n \rightarrow \infty} \inf \{ \mu (b) | b \in J^+(C, P(X)) \}$. Then $a$ has a subsequence $u$ that converges to some $v \in P(X)$. Now as $\cl (v) \cap C \neq \emptyset$, there is some $w \in IP(X)$ with $w \subset v $ and $\cl (w) \cap C \neq \emptyset$: Let $x \in \cl (v) $, $c \in C \cap \cl (I^-(x))  $ and $y: \N\rightarrow v$ with $y(n) \rightarrow_{n \rightarrow \infty} x$, let $c: x \leadsto$ past chr. curve, then by openness of the $I^-(y(n))$ it is easy to see that $c((0; t))) \in v$, thus with $w := I^-(c(0;t))$ we get the desired element of $IP(X)$. Now as $a$ was $\mu$-infimizing, a posteriori we see $v=w$, $w \in \cl (I^+(C, P(X))) $ and $w$ minimal: If $z <w$ then $d_1(\cl(z), \cl (w)) >0$, and $\cl (z) \setminus \cl(w)$ contains a ball, which is of positive measure, in contradiction to the assumption of an infimizing sequence. For the maximal case we correspondingly consider a $\mu$-supremizing sequence, a subsequence converging to some $V \in P(X) $, then each maximal indecomposable past subset $W$ of $V$ is a maximal element (here we cannot conclude $V=W$). \hfill \qed 

\begin{Theorem}[see \cite{oM14}]
\label{Respect}
The functor $F_+$ respects conformal future-compact extensions. More precisely, let $(M,g)$ be a globally hyperbolic spacetime and $E:(M,g) \rightarrow (N,h)$ be a conformal future-compact extension. Then the end-point map $\e_E$ is a homeomorphism from $(IP( M), \tau_+)$ to $ E(M) \cup \partial^+ (E(M),N)  \subset N$ taking $i_M(M)$ to $E(M)$. Its continuous inverse is the map $\overline{s}: p \mapsto E^{-1} (I_N^-(p))$. 
\end{Theorem}

\v {\bf Proof.} We want to show that $\overline{s}$ is a right and left inverse of $\e_E$. First we have to show that $\overline{s} \vert_{\partial^+ (E(M), N)}$ takes values in the IPs. Given a point $p$ in $\partial^+ (E(M), N)$, then there is a timelike future curve $c: [0;1] \rightarrow N$ from a point $q$ in $E(M)$ to $p$. Causal convexity of $\cl (E(M))$ implies that the curve is contained in $\cl (E(M))$ and that $k:= c \vert_{[0; 1)}$ even takes values in $E(M)$ and thus $ k= E \o \kappa$ for some $C^0$-inextendible timelike future curve in $M$. Then $E^{-1} (I_N ^-(p)) =  E^{-1} (I_N^- ( c)) = E^{-1} (I_N^- (E \o \kappa)) = I^-_M (k )$, where, in the last step, we use causal convexity. Therefore, indeed, for $q \in \partial^+ E(M)$, the set $E^{-1} (I^-(q) \cap E(M) ) $ is an IP in $M$. And $\overline{s}$ is a right inverse of $\e_E$ as $N$ is distinguishing, and a left inverse as $c$ generating the IP is a chr. chain.

The map $\overline{s}$ is continuous: First, the assignment $ p \mapsto I^- (p) $ is inner continuous in any regular chr. space (Theorem \ref{IC}) and is outer continuous in causally continuous chr. spaces, and g.h. manifolds are causally continuous (see \cite{BEE}. e.g.). Now let $a$ be a sequence in $N$ convergent w.r.t. the manifold topology. Then by i.o.-continuity of $I^-_N$ the sequence $k \mapsto I_N^-(a(k))$ and also the sequence $k \mapsto I_N^-(a(k) ) \cap E(M)$ converges in $\tau_{io}$ and thus also in $d_1$ (following Theorems \ref{HausdorffInnerOuter} and Theorem \ref{HB}). As $E$ is a homeomorphism onto its image, mapping compacta to compacta, the sequence $b:k \mapsto E^{-1} (I_N^-(a(n))$ converges in $\tau_{io}$ as well. Then Th. \ref{HausdorffInnerOuter} and Th. \ref{HB} imply convergence of $b$ in $d_1$ and Theorem \ref{TFAE} implies convergence of $b$ in $\tau_+$. For the converse direction, let $A_n \in IP(M)$, then each $A_n$ is a PIP $I^-(q_n)$ in $N$, so we want to show that convergence of $n \mapsto A_n$ implies convergence of $n \mapsto q_n $, but this is just the openness statement of Th. \ref{recovering}. \hfill \qed

\bigskip

A further positive feature of $\tau_+$ is that $\tau_+$ {\em inherits conformal standard-stationarity}, and mimicking the proofs in \cite{H0} and \cite{FHS1} one gets (see also the corresponding statement in \cite{CFH}) the following theorem. For its statement, recall that the {\bf Busemann function} $b_c: M \rightarrow \R \cup \{  \infty \}$ of an arclength-parametrized curve $c: I \rightarrow M$ is defined as $b_c (x) := \lim_{t \rightarrow \infty } (t - d(c(t), x)) $. The finite-valued Busemann functions are 1-Lipschitz, if we denote their set by $B(M)$ then the Busemann boundary is defined as $\partial_B M:= B(M) \setminus j(M)$ equipped with the compact-open topology, where for $p \in M$ we have $j(p) = - d(p, \cdot)$; for details see \cite{FHS2} and \cite{H}. The following theorem is from \cite{H}, we include a proof for the sake of self-containedness.

\begin{Theorem}[\cite{H}]
\label{StationaryRanders}
Let $(M,g)$ be a conformally standard static spacetime with standard slice $N$. Then ($\partial^+M, \tau_+) \setminus \{ M \} $ is homeomorphic to $ \partial_B N$.
\end{Theorem}

\V {\bf Proof.} As a Cauchy temporal function $t$ for the definition of the $f(t,A)$ (as in Theorem \ref{TFAE}) for the IPs $A$ we take the static temporal function whose level sets are the standard slices isometric to $N$. The Hausdorff distance in $\R \times N$ of any set $\{ (t,x) \vert t < b_c(x)  \}$ (which is $I^- (\hat{c} (\R))$ for $\hat{c}$ with $\hat{c} (t) := (t, c(t))$ being the unique null lift of $c$ starting at $t=0$) is after restriction to a compact subset equivalent to the $C^0$ (function space) topology on the Busemann functions, as the latter are $1$-Lipschitz. The rest is as in the cited articles (in particular Theorem 6 in \cite{H}). \hfill \qed

\medskip

\V Thus $IP(M)$ and $\partial^+ IP(M)$ are cones with tip $ M$, and the $\R$-action on $M$ mapping $(r, (t,x))$ to $(r+t,x)$ has a canonical continuous extension to $IP(M)$, which fixes the point $M$.   

\medskip

{\bf Remark:} A corresponding fact holds for the conformally standard stationary case, as in \cite{FHS2}.

\medskip

As an example, we consider a manifold $G$, which is a slight variant of Harris' unwrapped grapefruit-on-a-stick: $G = (\R^2, g_G:= (f^2 \o \pr_2) \cdot g_0)$ where $g_0 $ is the Euclidean metric and $f \in C^{\infty} ( \R , (0; \infty)) $ with $ f(\R \setminus [-2; 2] )  = 1$, $f(-x) = f(x) \forall x \in \R$, $f (x) \in (1;2) \forall x \in (-2; 2) \setminus [-1;1]$ and $ f(x) = 2 \forall x \in [-1;1]$, and $M:= (\R \times G , - dt^2 + g_G)$.

For every arc-length parametrized curve $c$ in $G$, we define a null curve $\hat{c}$ in $M$ by $\hat{c} (t) := (t, c(t))$. For every curve $c: I \rightarrow G$ whose Buseman function is finite, $c(t) / \vert c(t) \vert$ converges to some $s \in \SSS^1$. Conversely, if for two curves $c_1$, $c_2$, this limit is the same, the corresponding points $I^-(c_1), I^-(c_2) \in \partial^+ M$ are in the same $\R$-orbit. By distinction of the cases that $(\pr_2 \o c)^{-1} (0)$ is finite (in which case the curve stays in $\pr_2^{-1} ([-1;1])$) or infinite, we easily see that only curves $k$ in $M$ with $I^-(k) = M$ can $G$-project to 'oscillatory' curves in $G$, i.e., curves $\gamma$ intersecting $G \setminus B(0,r)$ for all $r >0$ in both semiplanes: $\forall r >0: \gamma (I) \cap  (G \setminus B(0,r)) \cap \R\times (\pm [0; \infty)) \neq \emptyset  $. Consequently, if for each 'asymptotic angle' $\lim_{t \rightarrow \infty} c(t) / \vert c(t) \vert \in \SSS^1$ we define $[c^\a] $ as the Busemann function of the standard Euclidean
 ray with angle $\a$ to the positive $x$-axis, $ d([c^\a], [c^\b]) >2$ whenever $\a \in (0; \pi) $ and $\b \in (\pi ; 2 \pi)$, and for $\a = 0$ or $\a = \pi$ there are two different classes of curves $[c_+ ^0]$ and $[c_-^0]$ with finite Busemann function, depending on whether the curve stays above or below $[-2;2]$). And also $d([c_+^0], [c_-^0]) \geq 2$, so the two classes are distant. Thus the Busemann boundary consists of 
 two cones over intervals $[0 ; \pi ]$ and $[ \pi; 2 \pi ]$ 
identified at their tips, parametrized by the asymptotic angle $\a$ 
(so at the $x$-axis at the values $0$ and $\pi$ each appear twice) and the $\R$-action.
Thus Theorem \ref{StationaryRanders} implies that $(\partial^+ M, \tau_+)$ consists of two cones over an interval identified at their tips (i.e., at $M$) as well. In contrast, $\tau_-$ is in this case non-Hausdorff (see \cite{FH}, whose Section 2 is emphatically recommended to the reader interested in the stationary case). Flores and Harris (in Sec. 2 of \cite{FH}, p. 1219 of the journal article, p.8 in the arXiv version) consider {\em (mutatis mutandis)} the sequence $\s: \N\rightarrow M$ defined by $\s (n) := (\frac{n}{2} + w, n,0)  $ for $w:= \int_{-2}^2 f(t) dt$, which is not a chr. chain, but still $\s (n) \in I^+(\s(0))  \forall n \in \N$. In $\tau_+$, the sequence $\s$ does not converge, showing that $(IP(X),\tau_+)$ is in general not future-compact. Even more, it shows that also $J^+(p) \cap \partial^+X  $ is in general noncompact, even if $X$ is g.h., as the $\s(n)$ can be made part of $\partial^+ X$ by considering $X:= M \setminus \bigcup_{k=1}^\infty J^+ (\s (k))$, which is an open chr. convex subset of $M$.

\bigskip

Here we see an interesting similarity to the case of a Riemannian manifold $(N,h)$. Its Gromov boundary $\partial_G N$ is defined as $(\cl (\{  [d_h (\cdot , x) ] \vert x \in N  \}, L_1 (N,d_h)/\R)) \setminus \{  d_h (\cdot , x) \vert x \in N \}$, where $[f]$ is the equivalence class of a function given by $f \sim \  g : \Leftrightarrow f-g$ constant, and $L_1 (N, d_h)$ is the space of $1$-Lipschitz functions on $N$ (cf. \cite{FHS2} for details). We call a Riemannian manifold $(N,h)$ {\bf (strongly) asymptotically $(k,l)$-Euclidean} iff there is a (compact) $K \subset M$ such that there is a diffeomorphism $D: N \setminus K \rightarrow U $ where $U = \R^n \setminus \cl (B(0, r_-))$ for some $r_->0$ s.t. in Euclidean coordinates $\partial_l ( g - g_0 )  \in \delta_{0s} g_0 + O(r^{-k-s})$ for each multi-index $s$ with $|s| \leq l$, where $g_0$ is the Euclidean metric. We call a Riemannian manifold $(N,h)$ {\bf (strongly) asymptotically Euclidean ((S)AE)} iff it is (strongly) asymptotically $(2,1)$-Euclidean. Each SAE manifold is complete. Examples for asymptotically Euclidean manifolds are the standard (Killing) Cauchy slices of (interior or exterior) Schwarzschild spacetimes $S_m$. More importantly, each $S_m$ is conformally equivalent to a standard static product with asymptotically Euclidean slices. They are not strongly asymptotically Euclidean due to the inner boundary.

\begin{Theorem}
	Let $(N,h)$ be SAE, $I \subset \R$ an interval and $M:= (I \times N, -dt^2 + h)$.
	\begin{enumerate}
		\item \( \partial_G N = \SSS^{n-1} = \partial_B N ,\)
		\item $(IP(M), \ll_{BS} , \tau_+)   = (IP(M), \ll_{BS}, \tau_-)  $, and $(IP(M), \tau_-) $ is homeomorphic to a (truncated if $I \neq \R$) cone over the Busemann completion $B(N)$ of $N$, the homeomorphism mapping the cone over $  \partial_B N  $ to $IP(M) \setminus i_M(M)$.
	\end{enumerate}

\end{Theorem}

 {\bf Proof.} For every proper Busemann function $b_c$ to a curve $c$, and more generally, for each sequence $A: \N\rightarrow N$ such that $d(t, \cdot )$ converges modulo additive constants, the map $ u_c: [0; \infty) \rightarrow \SSS^{n-1},  t \mapsto c(t) / r( c (t))  $ (rep. $ a: \N\rightarrow \mathbb{S}^{n-1}$, $a(n) := A(n) / r(A(n)) $) must converge, and $\lim_{t \rightarrow \infty} u_c (t)  = \lim_{t \rightarrow \infty } u_k (t)  $ implies $b_c = b_k$ for any two such curves (unlike in $\R \times G$!). 

Let $K>0$ with $||g - g_0|| \leq K/r$ on $r^{-1} ([r_-; \infty))$. Let $c$ be a Busemann curve. For $\phi, \psi \in \SSS^{n-1}   $ we fix two test points $ x_\phi = (r_0, \phi)$ and $x_\psi = (r_0, \psi )$ for $r_0>0$ s.t. for all $r>r_0 $ we have $|(g - g_0)_{\mu \nu}| < K/r  $. We wnt to estimate $| (d(x_{\phi} , (s, \phi )) - d(x_\psi, (s, \phi))) - (d(x_{\phi} , (t, \psi )) - d(x_\psi, (t, \psi)))|$ from below in terms of $\phi - \psi$ in order to bound the latter term (as $d (x_\phi ,y_n ) - d(x_\psi, y_n) \rightarrow_{n \rightarrow \infty} 0$ for a Gromov sequence $n \mapsto y_n$). Let w.l.o.g. $s \leq t $. For the envisaged estimate, we use four curves $c_s : u \mapsto (s, \gamma (u)), c_t : u \mapsto (t , \gamma (u)), k_\phi: u \mapsto  ((1-u) s + s u , \phi), k_\psi : u \mapsto ((1-u) s + u t, \psi) $. Then for the difference of the spherical parts we get $\ell (s) - \ell(r) \geq \a (s \sqrt{1-\frac{K}{s}} - r\sqrt{1 + \frac{K}{r}}) $. Let $f_\pm (x) := x \sqrt{1-K/x} $, then $f_\pm '(x) \rightarrow_{x \rightarrow \infty} 1$, so if $r_0$ is such that $ f'_pm (x) \in [9/10; 11/10]$ then $\ell(s) - \ell (r) \geq \a/2 \cdot (s-r)$. For the radial part we calculate $\sqrt{1+K/x} - \sqrt{1-K/x} = 2K (r \sqrt{1+K/x} + \sqrt{1-K/x}) \leq K/x$ for all $x \geq r_0$, thus $\ell (c_\phi) - \ell(c_\psi) \leq \int_r^2 K/x = K \ln (s) - K \ln (r()) $. That means: If the Gromov function for $a: \N\rightarrow N$ converges, then there is $\a \in \SSS^{n-1}$ such that $a (n) \in (n, \a) + O_{n \rightarrow \infty}(n^{-a}) $ for all $a <1$. Conversely, if we have this convergence for $a$ as above, then there is a curve $c$ with $b_c= g_a$. For example, take minimal curves $c_n$ from a fixed point to $a(n)$.  As the spherical line element is approximately proportional to $r$, we know that for the curve $k_\a: t \mapsto (t, \a)$ as a reference curve we get $ d (c(t) , k_\a (t) ) \leq H \cdot \sqrt{t} $ (take $a=1/2$ in the above estimate of the first part of the proof). Then we can use the estimate for $g'$ to get $\ell (c) - \ell(k) \leq \int \sqrt{r} \cdot r{-2} dr $, which converges in $\R$. 

For Item 2, note that if $\partial_B N$ is Hausdorff then $(IP(M), \tau_-)$ is a cone over $\partial_B N$ (\cite{FHS2}, Th. 1.2 Item 3A), thus Hausdorff, and if $\tau_-$ is Hausdorff, it agrees with $\tau_+$ (\cite{CFH}). \hfill \qed

\bigskip

Harris' example is very instructive as here $IP^+(X)= \cl (i_X(X), (P(X), \tau_+)) $ is strictly larger than $IP(X)$. As $\tau_-$ is coarser than $\tau_+$ we have $ \cl (i_X(X), (P(X), \tau_+)) $ strictly larger than $IP(X)$ as well. But each limit set of a sequence in $X$ w.r.t. the non-Hausdorff topology $\tau_-$ contains at least one element in $IP(X)$, thus $\partial^+X= \cl (i_X(X), (P(X), \tau_-)) \cap IP(X)$ is still future complete and sequentially future compact. In contrast, $(IP(X), \tau_+) $ is not future compact. This indicates that, using $\tau_+$, one should consider $IP^+(X)$ instead of $IP(X)$. As seen above, the two spaces coincide in the asymptotically flat case or in certain warped products over bounded intervals, cf. Section \ref{multiply-warped}.

\section{Application to multiply warped chronological spaces}
\label{multiply-warped}

We want to transfer statement and proof of \cite{H0}, Prop. 5.2, from $\tau_-$ to $\tau_+$, generalizing it at the same time from multiply warped spacetimes to multiply warped chronological spaces.

\bigskip

 For $a,b \in \{ - \infty \} \cup \R \cup \{ \infty \}  $, $(K_i, h_i)$ complete Riemannian manifolds and $f_i: (a;b) \rightarrow (0; \infty)$, the {\bf multiply-warped spacetime} $(a; b) \times_f \bigtimes_{j=1}^m K_j$ is defined as $(a;b) \times \bigtimes_{j=1}^m K_j$ with the metric $g:= -dt^2 + \sum_{j=1}^n (f_j \o \pr_0) \cdot (\pi_j^* h_j )$. One example of a multiply-warped spacetime is (interior) Schwarzschild spacetime with (after an obvious conformal transformation) $ I= (0; 2m) \ni r$, $K_1:= \R $, $K_2 := \mathbb{S}^2$, $f_1 (r) = (1-\frac{2m}{r})^2$, $f_2(r) = r^2 (\frac{2m}{r} -1)$ for all $r \in (0; 2m)$. The Grapefruit spacetime from Sec. \ref{MetricsOnIPs} is another example, as it is standard static. 

\bigskip

Let $I$ be a real interval. Let $(K_i, d_i)$ be complete length metric spaces and  $f_i \in C^1 (\Int (I) , (0; \infty)) $ for $i \in \N_m^*$, $X= I \times K, K:= \bigtimes_{i=1}^m K_i$, and let $ \pr_i : I \times K \rightarrow K_i $ be the standard projections, $\pr_0 = \pr_I : I \times K \rightarrow I$ and $\pr_K: I \times K \rightarrow K$. For $p \in X$, we write $p_i:= \pr_i (p)$. The topology of $X$ is the product topology $\tau_P$. There is a chronology $\ll_0$ defined by $(t,x) \ll_0 (s,y) $ if and only if there is a {\bf good path} from $(t,x)$ to $(s,y)$, which is a continuous path $c : [t;s] \rightarrow X$ from $(t,x)$ to $(s,y)$ with $c_0:= \pr_0 \o c \in C^1 $ with $ c_0(u) = u \ \forall u  $ and s.t. for all $r \in (t;s)$ we have a neighborhood $U (r)$ of $r $ s.t. each $\pr_i \o c \vert_{U(r)}$ is $k_{i,U}$-Lipschitz with respect to $f_i(r) d_i $, with $\sum_{i=1}^m k_{i,U}^2 < 1 (= (c'_0)^2 (u)  \ \forall u)$. A {\bf multiply warped chronological space over $I$} is a space $X = I \times K$ with the chr. relation $\ll_0$ via good paths described above. We denote this chr. space by $I \times^f \bigtimes_{j=1}^m K_j$ (mind the notation minimally differing from the one of a multiply-warped spacetime, only by upper resp. lower location of the symbol '$f$'). 

\bigskip

The next theorem shows that the basic relation between both notions: in a multi-warped spacetime, the future timelike curves are exactly the good paths. In this sense, multiply-warped spacetimes can be understood as multiply-warped chr. spaces.

\begin{Theorem}
\label{Persistence}
$\forall x,y \in I \times_f \bigtimes_{j=1}^m K_j : x \ll y \Leftrightarrow $ there is a good path $x \leadsto y$ in $I \times^f \bigtimes_{j=1}^m K_j$. 
\end{Theorem}

{\bf Proof:} A timelike curve $c: (u;v) \rightarrow I \times_f \bigtimes_{j=1}^m K_j$ can be reparametrized s.t. $c_0(t) = t$ and thus $c_0' =1$, then 

\[0 > g(c'(s), c'(s)) = -1 + \sum_{j=1}^{n} f_j(s) h_j (c_j'(s) , c_j'(s) ). \]

Any such $c$ is a good path: At each $r \in I$ we have 

\[0 > g(c'(r) , c'(r))  = -1 + \sum_{j=1}^m f_j(r) h_j (c_j' (r) , c_j '(r)). \]

By continuity, there is $\e >0$ s.t. 

$$ \sum_{j=1}^m \max \{ f_j (r) g(c_j'(u) , c_j' (u) ) \vert u \in [r - \e ; r + \e ] \}  <1 ,$$

so $c$ is good. Conversely, let a good path $c$ from $x$ to $y$ be given, then for all $r \in (a;b)$ we want to show that there is $\e >0$ with $c(r-u) \ll c(r ) \ll c(r+u) $ for all $u \in (0; \e]$. For $U (r)$ as in the definition of goodness of a curve, let $\kappa:= c\vert_{U(r)}$. Then the $\kappa_i$ are $k_i$-Lipschitz w.r.t. $f_i (r) d_i $ for all $i \in \N_m^*$ and $ \rho = \sum k_i ^2 <1$. We find $ \d >0$ with $(1 + \d) \rho <1$ and we define $\e >0$ s.t. $f_i (c(x)) < (1 + \d) f_i (c(r)) \ \forall u \in [r - \e ; r + \e]$. Then we can construct the required timelike curve $\g$ between $c(r) $ and $c(r+u)$ as a $d_i$-geodesic curve in each factor $K_i$ From $h_i (\g_i ' (x) , \g_i ' (x)) < k_i^2$ we conclude

$$ g(\g_K'(x) , \g_K'(x)  ) = \sum_{i=1}^m f_i (x) h_i (\g_i' (x) , \g'_i (x) ) < (1 + \d) \sum_{i=1}^m k_i <1 , $$

which shows that $\g$ is timelike. \hfill \qed

\begin{Theorem}
Let $ X  = I \times^f \bigtimes_{j=1}^m K_j$, $I= (a;b)$ with $b < \infty$\footnote{This has been choosen for the sake of simplicity of the presentation and can be arranged by a reparametrization without disturbing the other properties.}. Assume that for all $i \in \N^*_m$, $(K_i, d_i) $ is a complete length space (i.e. each Cauchy sequence in $K_i$ converges), and
 (*)  $\int_c^b f_i^{-1/2} (x) dx < \infty$ for some (hence any) $c \in (a;b)$. 
Then $(IP(X) = IP^+(X), \ll_{BS}, \tau_+)$ is chronologically homeomorphic to $ (\ov{X} :=  (a; b]  \times K_1 \times ... \times K_n, \ll_0)$.
\end{Theorem}

\V{\bf Remark.} Schwarzschild spacetime satisfies the additional finiteness condition (*) on the $f_j$, the Grapefruit Spacetime does not because of the infinite range of the time parameter. 

\bigskip

\V {\bf Proof}. The map $L: x \mapsto I^-(x)$ is a chr. isomorphism from $ ((a; b] \times K, \ll_0) $ to $(IP(X) , \ll_{BS}) $: The map $i_X(p) \mapsto p, I^-(c(\N)) \mapsto \lim_{n \rightarrow \infty}^{\ov{X}} $ is inverse to $L$, Eq. (*) implies that for every chr. chain $c$ in $X$, $\pr _K \o c$ converges to some $k_c \in K$, thus to each $c$ without endpoint in $X$ we can assign the point $(b, k_c)$, and as $i_X$ is chronological, we only have to check whether $y \gg_0 i_X(z) $ for all $y \in \partial^+ IP(X)$ and all $z \in y$, but this is true as every chr. chain in $X$ can be interpolated and equipped with an endpoint as above to become a good path in $\ov{X}$.

Both sides are metrizable. Now it is an easy exercise to see that convergence of the closures of past sets in $X$ is coordinate-wise, i.e. for a sequence $A: \N\rightarrow P(X)$ we define 

\[B:= \limsup_{n \in N} (\sup \{ \pr_0 (x) | x \in A_n \})\]

and distinguish the cases $B <b$ and $B=b$ and show that in each case the limsup is a true limit, and that, for $ i_X(X) \ni A(n) = I^-(p_n) $, the projection to the other factors have to converge as well, showing also that $IP^+(X) = IP(X)$. \hfill \qed

 \newpage

\section{Causal perspective, solomonic conclusion, open questions}
\label{posets}

Finally, one can look at the causal completion from the perspective of the causal relation $\leq$ instead of the chr. relation $\ll$. One obvious advantage of this procedure is that we can require semi-fullness (of $\leq$) right from the beginning and still include causal diamonds $J^+(p) \cap J^-(q)$.

Recall that there is a functorial way to induce a causal (an order) relation from a chr. one, a map $\a : X \times X \rightarrow X \times X $ defined by 
 \( x \ \a (\ll) \   y : \Leftrightarrow I^+(y) \subset I^+(x) \land I^-(x) \subset I^-(y) \forall x,y \in X \).

{\em Enriched by adding a causal relation $\leq = \a$, the functor $IP_+$ maps every causally continuous well-behaved chr. space $X$ to a causally simple chr. space $IP(X)$}: 

Let $p \in X$ and $a: \N \rightarrow J_{IP(X)}^+(p)$, i.e., $ a(n) \supset p \forall n \in \N $ and $a(n) \rightarrow_{n \rightarrow \infty}^{\tau_+} v $, i.e., $I_X^-(\liminf a) = v = I_X^-(\limsup a)$. Then $ p = I^-_X (p) \subset I_X^- (\cap_{j \in \N } a(j)) \subset I_X^- (\liminf a) = v$. The case $ a: \N \rightarrow J_{IP(X)}^- (p) $ is proven analogously, here the last line is $p= I_x^- (p) \supset I_X^- (\bigcup_{j \in \N} a(j)) \supset I_X^- (\limsup a) = v$. 

This corresponds to the fact mentioned above that $\a$ introduces spurious causal relations in $i_X(X)$. 

Fortunately, we can refine $\a$ to a functor $\tilde{\a}$ from chr. to ordered spaces $(X, \ll) \mapsto (X,\tilde{\a} (\ll))$ with $\a (\ll) \supset \tilde{\a} (\ll)$ that avoids this undesired phenomenon: We define, for any real interval $I = [a;b]$ and any chr. space $(X, \ll)$ a curve $c: I \rightarrow X $ to be {\bf chr. coherent} iff 

$\forall t \in (a;b]: I^-(c(t)) = \bigcup_{s \in [a;t)} I^-(c(s)) \ $ and $\forall t \in [a;b) : I^-(c(t)) = I^- \big ( \bigcap_{s \in (t; b]} I^-(c(s))  \big) \  $, and

$$x \tilde{\a}(\ll) y : \Leftrightarrow  x \ll y \lor \big( \exists c: [0;1] \rightarrow X :  c(0) = x  \land c(1) = y \land c {\rm \ chronologically \ coherent} \big)  . $$    

The definition of $\leq = \tilde{\a} (\ll)$ implies $\ll \subset \tilde{\a}(\ll) \subset \a (\ll) $ and the {\bf push-up property}
$x \leq y \ll z \lor x \ll y \leq z \Rightarrow x \ll z \forall x,y,z \in X$
(valid even for $\a(\ll) $ and in every spacetime, see e.g. \cite{eMmS}).

As both $\tau_-$ and $\tau_+$ are chr. dense, the push-up property entails $\Int (J^\pm(A)) = I^\pm (A)$ for each $A \subset X$ for both topologies.

An advantage of $\tilde{\a}$ over $\a$ is the absence of spurious causal relations for causally continuous spacetimes even if they are not causally simple:

\begin{theorem}
	Let $(X, g)$ be a causally continuous spacetime with $\ll \subset M \times M$ resp. $\leq \subset M \times M$ being the usual chr. resp. causal relation defined by connectability along timelike resp. causal curves. Then $ \tilde{\a} (\ll) = \leq$. 
\end{theorem}

\V{\bf Proof.} The inclusion of the right-hand side in the left-hand side follows from the fact that each $\leq$-causal curve is an $\a (\ll)$-chr. coherent curve. For the other inclusion, let $p, q \in M$ with $p \tilde{\a} (\ll) q$. Then there is an $\a(\ll)$-chr. coherent curve $c: I \rightarrow M,  p \leadsto q$. Now $c$ is continuous: Let a sequence $ t : \rightarrow I $ with $t (n) \rightarrow_{n \rightarrow \infty} t_\infty \in I$. Then any subsequence $t \o j$ has a monotonous subsequence $s:= t \o j \o k $, such that $ I^- (c \o s) $ converges to the union or intersection of the respective subsets equalling $I^- (c(t_\infty))$, and in the manifold topology (which equals $\tau_+$) any sequence of points whose pasts converge converges. Then $c$ is also a (continuous) causal curve, as every spacetime is locally causally simple, which allows to use \cite{eMmS} Th. 3.69. \hfill \qed

\bigskip

A {\bf causal set}\footnote{The author is well-aware of the fact that this notion deviates from the one used in causal set theory by not requiring local finiteness. Nonetheless, the analogy to chr. sets is too striking to renounce using temporarily this notion in the context of our article.} is a (partially) ordered set $(X, \leq) $ that is semi-full and causally separable (defined like chr. separable, only replacing $\ll$ with $< = \leq \setminus \Delta$) and nowhere total, i.e., $\forall x \in X : X \neq J(x)$. 

$X$ is called {\bf future} resp {\bf past causally simple} iff $J^\pm (p) $ is closed for every $p \in X$, or equivalently, that $J^+ \subset X \times X $ is closed: The nontrivial implication '$\Rightarrow$' follows from $J^+ = \cl (I^+) $, whose nontrivial inclusion '$\supset$' follows from $(x,y) \in \cl (I^+) \Rightarrow y \in \cl (I^+(x)) = J^+(x)$.

Analogously to the above, a {\bf causal space} is a causal set equipped with a topology. A preregular causally continuous, locally causally simple path-generated causal space $ (X, \leq , \tau) $ is called {\bf well-behaved} iff $(X, \tau)$ is locally compact, $\sigma$-compact and locally arcwise connected.

Conversely to the functors $\a$ and $\tilde{\a}$, there is a functorial way $\beta$ to define a chr. relation from a causal one. Following a suggestion of Miguzzi and S\'anchez \cite{eMmS} (Def. 2.22), $ \b (\leq)$ defined by 

$$(x,y) \in \b (\leq) : \Leftrightarrow \big( x \leq y \land ( \exists u,v \in X:  x < u< v < y \land J^+(u) \cap J^-(v) {\rm \ not \ totally \ ordered} ) \big)$$

is a chr. relation which encodes the fact from Lorentzian geometry that a null geodesic ceases to be maximal after the first cut point (thus $\beta$ recovers the usual chr. relation in each causal spacetime, cf. \cite{eMmS}(Th.3.9)) and manifestly has the {\em push-up property}. Under which conditions on chr./causal sets $\beta $ is inverse to $\tilde{\a}$ is an interesting question that cannot be pursed further here. The warning example of $\R^{1,2} \setminus I(0)$ shows that semi-fullness of $\b(\leq)$ is not automatically satisfied. The brute-force way here would be to just assume a partial order relation $\leq$ on $X$ such that $\b(\leq) $ and $\tau_+$ make $X$ a well-behaved chr. space. But already with much weaker hypotheses we can show that for a well-behaved causally continuous causal space $X$, the causal space $(ICP(X), \subset, \tau_+ )$ of indecomposable causal past subsets is well-behaved and causally continuous by the proof of Theorem \ref{recovering} {\em mutatis mutandis}. Again, all central statements work also in the subset $CP(X) $ of causal past subsets instead of the indecomposable ones $ICP(X)$.

\bigskip

The causal perspective to induce $\ll$ from the primordial object $\leq$ allows for a reinterpretation of $IP(X)$ as the set of filters on the poset $(X, \leq) $, and for a simple one-line description for a future completion as a causal-chr. space with a metric $\delta$, where $d_1$ is the above pointed Hausdorff metric: 

$$ (IP(X), \leq := \subset, \ll := \beta (\leq), D := \cl^* (d_1)) .$$

 $IP(X)$ is $J^-$-distinguishing, as $x = \bigcup \{ J^-(x) \vert  x \in IP(X) \} $. The push-up property entails  $\bigcap \{ I^+(r) \vert r \in R \} = \bigcap \{ J^+(r) \vert r \in R \}$ for every past subset $R$ of $X$, and consequently we can and will define $\tau_+$ in our new context using $J^{+ \cap} $ instead of $I^{+ \cap}$. The following theorem shows that the diagram

\begin{equation}
	\label{CommDiag}
	\xymatrix{
		(X, \ll) \ar@{|->}[r]^{} \ar@{|->}[d]_{\tilde{\a}} & (IP(X), \ll_{BS}) \ar@{|->}[d]_{\tilde{\alpha}} \\
		(X, \leq) \ar@{|->}[r]^{} & (IP(X), \subset) 
	}
\end{equation}

is commutative, but one should keep in mind that the construction of $IP(X)$ also in the lower line involves the chr. relation, which we will see can defined as $\beta(\leq)$ there. Alternatively, we can define causal past sets $CP(X)$ and define indecomposability analogously to the chr. case.

\begin{Theorem}
	For every chr. set $X$, we have $ A \ \a(\ll_{IP(X)}) \ B \Leftrightarrow A \subset B $.
\end{Theorem}

{\bf Proof.} For the implication from the right to the left, assume $U \ll_{BS} A$, then $\emptyset \neq I^{+ \cap}_X (U) \cap A \subset I^{+ \cap}_X (U) \cap B$, thus $U \ll_{BS} B$, so $I_{IP(X)}^-(A) \subset I_{IP(X)}^-(B)$. But also $I^+_{IP(X)} (B) \subset I^+_{IP(X)} (A)$: For $U \in I^+_{IP(X)} (B)$ we have $\emptyset \neq I_X^{+ \cap} (B) \subset  I_X^{+ \cap} (A) $ due to the monotonicity of $I_X^{+ \cap}$. For the reverse implication, let $A \a(\ll_{BS}) B$, i.e. $I^-_{IP(X)}(A) \subset I^-_{IP(X)} (B)$, so for every $U \in IP(X)$ we have $I_X^{+ \cap } (U) \cap A \neq \emptyset \Rightarrow I_X^{+ \cap } (U) \cap B \neq \emptyset$. Now let $a \in A$. Then $U:= I^-(a) \in IP(X)$ and $a \in I^{+ \cap }_X(U) \cap A \neq \emptyset  $, so $ I^{+ \cap }_X (U) \cap B \neq \emptyset$ and $ a \in B$.   \hfill \qed

\bigskip

\V We define a chr. space $X$ to be {\bf globally hyperbolic} or {\bf g.h.} iff $J^+(K) \cap J^-(C)$ is compact for any two compact sets $C,K$ in $X$.
In the globally hyperbolic case, the future boundary is achronal:

\begin{Theorem}
\begin{enumerate}
\item Let $(X, \ll, \tau_+) $ be a well-behaved past- and future-full g.h. chr. space. Let $p \in IP(X) \setminus i_X(X)$. Then $J^+(p) \cap i_X(X) = \emptyset = I^+(p)$ for $(\ll_{IP(X)}, \leq_{IP(X)}) = (\ll_{BS}, \a(\ll_{BS}) = \subset)  $.
\item Let $(X, \leq)$ be a causal set s.t. $(X, \leq, \beta(\leq) , \tau_+)$ is a well-behaved g.h. chr. space. Let $p \in IP(X) \setminus i_X(X)$. Then $J^+(p) \cap i_X(X) = \emptyset $, for $(\ll_{BS}, \leq_{P(X)}) = (\beta(\subset), \subset) $.
\end{enumerate}
\end{Theorem}

\V {\bf Proof.} From Sec. 2 we know that $i_X(X)$ is past in $IP(X)$, which, together with future-fullness of $X$ and the push-up property, settles the first equality in the two items. From Th. \ref{recovering} we know that, in the setting of Item 1, $I_{IP(X)}^+$ is open and $c(n) \rightarrow c_\infty = \bigcup_{n \in \N} c(n)$ for a chr. chain $c$ in $IP(X)$. This is, if $q \in I^+(p) $ for $ p \in IP(X) \setminus i_X(X)$, let $c$ be a chr. chain with $c(n) \rightarrow_{n \rightarrow \infty} q$, then there is $m \in \N$ s.t. $c(m) \in I^+(p)$, contradicting $J^+(p) \cap i_X(X) = \emptyset$. \hfill \qed  

\bigskip

The next theorem shows that the causal relation defined by $\a $ or $\tilde{\a}$ does not deviate too much from the chr. one even in chr. spaces: 

\begin{Theorem}
	\label{JandI}
	
		Let $(Y, \ll, \tau)$ be a chr. dense chr. space and $\leq = \a (\ll) $ or $\leq = \tilde{\a} (\ll)$ then we have $ J^\pm(p) \subset \cl (I^\pm(p)) \forall p \in X$.

\end{Theorem}

{\bf Proof:} Directly from the push-up property. \hfill \qed

\begin{Theorem}
	\label{Ladder-piece}
	\label{CSCC}
	Let $(X, \ll, \tau)$ be a chr. dense regular chr. space. Then \\ $x \a (\ll) y \Rightarrow i_X(x)  \a (\ll_{BS}) i_X (y) $. Moreover, if $(X, \ll, \tau)$ is $\a$-causally simple then
	\begin{enumerate}
		\item $x \a (\ll_X) y \Leftrightarrow i_X(x)  \a (\ll_{BS}) i_X (y) $.
		\item $J^\pm (C) = \cl (I^\pm (C))$ for all compact subsets $C$ of $X$, in particular $J^\pm (p) = \cl (I^\pm (p)) \forall p \in X$.  
		\item $ (X, \ll, \tau)$ is causally continuous.
		
	\end{enumerate}
	
\end{Theorem}

\V {\bf Proof.} The first assertion follows from $x \leq y \Rightarrow I^-(x) \subset I^-(y) \Rightarrow i_X(x) \leq i_X(y)$: For $A,B \in P(X)$ we want to show $A \subset B \Rightarrow A \a (\ll_{BS}) B$. And indeed, $A \subset B$ is equivalent to 

\[ \forall C \subset X : (\exists x \in C: I^-(x) \supset B \Rightarrow \exists x \in C: I^-(x) \supset A)  \ \land \ ( \exists x \in A: I^-(x) \supset C \Rightarrow \exists x \in B: I^-(x) \supset C) , \] 

the direction from left to right being obvious, for the other one choose $y \in A$ and $C:= I^-(y)$. 
 
The second item applied to points follows from the the first assertion of Th. \ref{path-generated}. 
For the second item applied to general compact sets, let $y \in \cl (J^+(C))$, then there is a sequence $a: \N \rightarrow X$ and a sequence $b: \N \rightarrow C$ with $ a(n) \rightarrow_{n \rightarrow \infty} y $ and $b(n) \leq a(n) \ \forall n \in \N $. Then a subsequence $b \o j $ of $b$ converges to some $x \in C$. As causal simplicity implies that $J^+ \subset X \times X$ is closed and the (in $X \times X$) convergent sequence $(a, b \o j) $ takes values in $J^+$ we can conclude $ y \in J^+(x) \subset J^+(C)$. So $J^+(C)$ is closed, and obviously $J^+(C) \subset \cl (I^+(C)) $, which implies the assertion.

For the third item, $\cl (I^\pm (x)) = J^\pm (x)  \ \forall x \in X$ implies $\forall p,q \in X: q \in \cl (I^\pm (p) ) \Rightarrow p \in \cl (I^\mp (q))$ (*). Now assume that $I^-$ is {\em not} outer continuous, then there is $p \in X$ and a compact $C \subset X$ such that $C \cap \cl (I^- (p)) = \emptyset $ (**) yet there is a sequence $a$ with $a (n) \rightarrow_{n \rightarrow \infty} p$ such that $C \cap \cl (I^-(a(n))) \neq \emptyset \ \forall n \in \N$. Choose $c(n) \in C \cap \cl (I^-(a(n)))$. Then a subsequence $c \o j $ converges to some $c_\infty \in C$. For all $y \ll c_\infty$, openness of $I^+(y)$ implies that $y \ll a(n)$ for $n$ sufficiently large, thus $p \in \cl (I^+(y))$. This implies with (*) that $y \in \cl (I^-(p))) = J^-(p)$. Now let $z (n) \rightarrow_{n \rightarrow} c_\infty$ with $z (n) \ll c_\infty \ \forall n \in \N$, then $z(n) \in J^-(p) \ \forall n \in \N$, which together with closedness of $J^-(p) $ implies that $c_\infty \in J^-(p) $. As $c_\infty \in C$, we have $J^-(p) \cap C \neq \emptyset$, contradicting (**).

For the missing implication in the first item, consider $i_X (x) \leq i_X(y) \Rightarrow_{Eq. \ref{CommDiag}} I_X^-(x) \subset I_X^- (y) $, which implies by chr. denseness $x \in \cl (I^-(y)) = J^-(y)$.  \hfill \qed

\medskip

To transfer the results to $\tilde{\a}$, note that a causally coherent path maps to a causally coherent path under $i_X$, and that Th. \ref{JandI} implies that every $\tilde{\a}$-causally simple chr. space is $\a$-causally simple.

\V Now we want to define strong causality for chr. spaces. This point is a bit subtle, as the embedding of strong continuity on the causal ladder relies on the non-imprisonment property for causal curves, which in turn is proven via local coordinates. The same holds for the proof that any strongly causal spacetime  carries the Alexandrov topology (for both see \cite{Penrose:Techniques}). 

A chr. space $X$ is {\bf almost strongly causal} iff every neighborhood of any $x \in X$ contains a causally convex subneighborhood of $x$.  A causal chain $c$ in an almost strongly causal chr. space $X$ converges to $p \in X$ if and only if $I^-(c) = I^-(p)$.
$X$ called {\bf strongly causal} iff every neighborhood $U$ of any $x \in X$ has a subneighborhood of the form $J^+(K) \cap J^-(L)$ for compact subsets $ K,L \subset X$ (we can even choose $K, L \subset U$ if $X$ is locally compact, locally path-connected and path-generated). If $X$ is causally simple, by intersection with $J^\pm(x)$ we can ensure that $K \subset J^-(x)$, $L \subset J^+(x)$ if $J^+(x), J^-(x) \neq \emptyset$. Finally, $X$ is called {\bf Alexandrov} iff every neighborhood of any point $x \in X$ has a subneighborhood of $x$ of the form $J^+(y) \cap J^- (z) $ for some $y,z \in X$, i.o.w., iff its topology is the Alexandrov topology. The three notions are different: The weakest one, almost strong causality, is equivalent to requiring that the subneighborhood can be chosen to be the intersection of a future and a past subset (write a convex subneighborhood $V$ as $I^+(V) \cap I^-(V)$). Secondly, $X_0:= \R^{1,2} \setminus I^+(0)$ wih the induced chronology and topology is strongly causal (even causally simple) but not Alexandrov. For an example of an almost strongly causal yet non-strongly causal chr. space, replace spatial $\R^2$ in $X_0$ with the infinite-dimensional separable Hilbert space $\ell_2$, losing local compactness.

It is easy to see that for $X$ well-behaved, we have $\Int (J^\pm (A) ) = I^\pm (A) $ for all $A \subset X$ (where either $X$ is a chr. space and $\leq := \tilde{\a} (\ll)$ or $X$ is a causal space and $\ll := \beta (\leq)$ or $\ll := \gamma(\leq)$).

If $X$ is causally continuous, then $X$ is strongly causal, and $\underline{X}$ is Alexandrov: Let $x \in U_2 \subset X$ be open. We define $U_3$ and $L$ as in the proof of Th. \ref{CSimpliesSC} while replacing $CB(x,r)$ with $A(r):= J^+(c(-t)) \cap J^-(c(t))$ for a $C^0$ future timelike curve $c$ with $c(0)=x$, then we argue as in the proof of the last item of that theorem, as $ A(0) = \{ x \}$.

So if $(X, \leq )$ is a causal set with $\a \o \gamma (\leq)  = \leq $ and such that $(X, \leq, \gamma (\leq), \tau)$ is a well-behaved past-full g.h. chr. set, then for $ICP(X)$ we have $(ICP(X), \subset, \gamma(\subset) ) = (IP(X), \subset, \ll_ {BS} ) $ where $ \ll_{BS}$ is defined via the chr. relation $\gamma(\leq) $ on $X$: We have the commuting diagram

\begin{equation}
	\label{CommDiag}
	\xymatrix{
		(X, \leq) \ar@{|->}[r]^{} \ar@{|->}[d]_{\gamma} & (ICP(X), \subset) \ar@{|->}[d]_{\gamma} \\
		(X, \ll = \g(\leq)) \ar@{|->}[r]^{} & (IP(X), \ll_{BS}) 
	}
\end{equation}

where we identify $A \in ICP(X)$ with $I^-(A) \in IP(X)$.  

\medskip

The assignment of the topologies  $\tau_\beta := \tau_+ \o \beta$ and $\tau_\gamma := \tau_+ \o \gamma$ constitutes functors from the category POS of partially ordered sets to the category TOP of topological spaces, and, as opposed to the Alexandrov topology, yielding the desired topology on causal spaces with future boundary.

\medskip

After so many facts about $\tau_-$ and $\tau_+$, our question remains which topology is more advantageous. Unsurprisingly, the author is convinced that the choice between $\tau_-$ and $\tau_+$ has revealed itself more and more as a pure matter of taste. Those who enjoy applicability to the wide class of strongly causal spacetimes, future compactness and compactifications to both time directions (keep in mind that $M \cup \partial^+M$ alone is noncompact for both topologies!) might prefer $\tau_-$, whereas those attracted by metrizability and interested in the future completion as a quick practical tool are likely to tend to $\tau_+$, in particular if one widens the scope to include all of $P(X)$. Theorem \ref{ConvergencePlusMinus} shows that the two topologies on $IP(X)$ carry basically the same information (whereas on $P(X)$ or even $\partial^+ P(X)$  the topology $\tau_+$ carries slightly but strictly more information, so whereas $(IP(X), \tau_-)$ is a good {\em minimal} future completion, maybe $(P(X), \tau_+)$ is a good larger future completion with some additional advantages). An interesting question is whether there is a holographic principle at work allowing to reconstruct the chr. structure of $X$ from $(\partial^+ X, \ll_{BS}, \tau_-)$ or $(\partial^+ X, \ll_{BS}, \tau_+)$.
 
In a recent work \cite{oM-Finiteness}, the author proved a Lorentzian finiteness theorem, defining a functor from a Lorentzian to a Riemannian category and using Cheeger finiteness, which in turn uses Gromov compactness of the space of compact metric spaces. It is quite obvious that such a functor is restricted to timelike compact Lorentzian manifolds, e.g. Cauchy slabs. But there is an refined approach to Gromov compactness: metric measure spaces, of which Berger writes in his book \cite{mB} that they are, in his opinion {\em the geometry of the future}. The causal perspective shows that they could have an appropriate Lorentzian analogue: {\em (partially) ordered measure spaces}, forming a category $PM$, whose objects comprise Lorentzian spacetimes but also possible degeneracies.  There is an obvious injective functor from the category of causal Lorentzian manifolds into $PM$: The order allows to reconstruct the conformal structure and together with the volume form one obtains the entire Lorentzian metric. To obtain finiteness theorems, it is fundamental to know compactness of appropriate subsets of ${\rm Obj}(PM)$, subject of a row of forthcoming articles starting with \cite{oM-Finiteness}, \cite{oM-LorGH}.

\bigskip

%{\bf Statement:} For this article, there is no conflict of interests.

\newpage

\begin{small}

\end{small}

\end{document}